%% Major changes done on March 30, 1996
%% TWO THEOREMS, NEWEST VERSION
%% Author: ALEKSANDAR IVI\'C
%% Input: plain TeX
%% Text typed Beograd December 1995 - Limoges January 1996 
%% - Beograd February-April 1996
%% Presented at the Bordeaux Conference in honour of the
%% Prime Number Theorem Centenary, Bordeaux, January 26, 1996
%% This is now the final version
%% for the Kyoto Conference on Analytic Number Theory, May, 1996.

%% Revision for ArXiv Mathematics in November 2003

\vglue 2cm
\nopagenumbers

\pageno=1

  \def\rightheadline{{\hfil{\tenrm
  On some reasons for doubting the Riemann hypothesis
  }\hfil\tenrm\folio}}

  \def\leftheadline{{\tenrm\folio\hfil{\tenrm
  Aleksandar Ivi\'c }\hfil}}
  \def\emptyheadline{\hfil}
  \headline{\ifnum\pageno=1 \emptyheadline\else
  \ifodd\pageno \rightheadline \else \leftheadline\fi\fi}
  \font\aa=cmss12
  \font\bb=cmcsc10
  \font\cc=cmcsc8
  \font\dd=cmr12
  \font\ee=cmtt8
\def\DJ{\leavevmode\setbox0=\hbox{D}\kern0pt\rlap
{\kern.04em\raise.188\ht0\hbox{-}}D}

\def\txt#1{{\textstyle{#1}}}
\baselineskip=13pt
\def\hf{{\textstyle{1\over2}}}

\def\d{{\,\rm d}}

\def\G{\Gamma}

\def\={\;=\;}

\def\zt{\zeta(\hf+it)}

\def\no{\noindent}  
  
\def\z{\zeta}

  \centerline{\aa ON SOME REASONS FOR DOUBTING THE
  RIEMANN HYPOTHESIS}
  \bigskip \medskip
  \centerline{\bb Aleksandar Ivi\'c}
  \bigskip \medskip
  {\cc  Abstract}. {\it Several arguments against the truth of the Riemann
  hypothesis are extensively discussed. These include the Lehmer
  phenomenon, the Davenport--Heilbronn zeta-function, large  and
mean  values  of} $ |\zt|$ {\it on the critical line,
and zeros of a class of
convolution functions closely related to  $ \zeta (\hf + it) $.
 The first two topics are classical, and
the remaining ones are connected with the author's recent research.
By means of convolution functions a conditional disproof of the Riemann
hypothesis is given.}
  \bigskip \smallskip
  {\dd 0. Foreword (``Audiatur et altera pars")}
  \bigskip\smallskip
This is the unabridged version of the work that was presented
at the Bordeaux Conference in honour of the
Prime Number Theorem Centenary, Bordeaux, January 26, 1996
and later during the 39th Taniguchi International Symposium
on Mathematics ``{\it Analytic Number Theory}", May 13-17, 1996
in Kyoto and  its forum, May 20-24, 1996. The abridged printed version,
with a somewhat different title, is [62]. 
The multiplicities of zeros are treated in [64].
A plausible conjecture for the coefficients of the main term in the
asymptotic formula for the $2k$-th moment of $|\zt|$ (see (4.1)--(4.2))
is given in [67].

\medskip
In the years that have passed after the writing of the first
version of this paper, it appears that
the subject of the Riemann Hypothesis has
only gained in interest and importance. This seem particularly
true in view  of the Clay
Mathematical Institute prize of one million dollars for the proof
of the Riemann Hypothesis, which is called as one of the mathematical
``Problems of the Millenium". A comprehensive account is to be found
in E. Bombieri's paper [65]. It is the author's belief that the present 
work can still be of interest, especially since the Riemann Hypothesis
may be still very far from being settled. Inasmuch the Riemann Hypothesis
is commonly believed to be true, and for several valid reasons,
I feel that the arguments that disfavour it should also be pointed
out.

\medskip
 One of the reasons that
the original work had to be shortened and revised before being published
is the remark that ``The Riemann hypothesis is in the process of being
proved" by powerful methods from Random matrix theory (see e.g.,
B. Conrey's survey article [66]).
Random matrix theory has undisputably found its
place in the theory of $\z(s)$ and allied functions (op. cit. [66], [67]).
However, almost ten years have passed since its advent, but the
Riemann hypothesis seems as distant now as it was then.

\medskip

  \bigskip \smallskip
  {\dd 1. Introduction}
  \bigskip\smallskip
  A central place in Analytic number theory is occupied by the Riemann
  zeta-function $\zeta(s)$, defined for $\Re{\rm e}\,s > 1$ by
  $$
  \zeta(s) = \sum_{n=1}^\infty n^{-s} =
  \prod_{p \,{\rm prime}} (1 - p^{-s})^{-1}, \leqno{(1.1)}
  $$
  \noindent
  and otherwise by analytic continuation. It admits meromorphic
  continuation to the whole complex plane, its only singularity
  being the simple pole $s = 1$ with residue 1. For general
  information on $\zeta(s)$ the reader is referred to the
  monographs [7], [16], and [61]. From the functional equation
  $$
  \zeta(s) = \chi(s)\zeta(1 - s), \;\chi(s) = 2^s\pi^{s-1}\sin\bigl(
  {{\pi s}\over 2}\bigr)\Gamma(1 - s), \leqno{(1.2)}
  $$
  \noindent
  which is valid for any complex $s$, it follows that $\zeta(s)$
 has zeros at $s = -2, -4, \ldots\,$ . These zeros are traditionally
 called the ``trivial" zeros of $\zeta(s)$, to distinguish them
 from the complex zeros of $\zeta(s)$, of which the smallest ones
 (in absolute value) are ${1 \over 2} \pm 14.134725\ldots i$.
 It is well-known that all complex zeros of  $\zeta(s)$  lie
 in the so-called ``critical strip" $0 < \sigma = \Re{\rm e}\, s < 1$, and
 if $N(T)$ denotes the number of zeros $\rho = \beta + i\gamma$
($\beta, \gamma$ real) of  $\zeta(s)$ for which $0 < \gamma \le T$,
then
$$
N(T) = {T\over {2\pi}}\log\bigl({T\over {2\pi}}\bigr) - {T\over {2\pi}}
+ {7\over 8} + S(T) + O\left({1\over T}\right) \leqno{(1.3)}
$$
\noindent with
$$
S(T) = {1\over \pi}\arg \zeta(\hf + iT) = O(\log T). \leqno{(1.4)}
$$
This is the so-called Riemann--von Mangoldt formula.{\it The Riemann
hypothesis} (henceforth RH for short) is the conjecture,
stated by B. Riemann in his epoch-making memoir [52], that
{\it very likely all complex zeros of  $\zeta(s)$
have real parts equal to} 1/2. For this reason the line
$\sigma = 1/2$ is called the ``critical line" in the theory of  $\zeta(s)$.
Notice that Riemann was rather cautious in formulating the RH, and
that he used the wording ``very likely" (``sehr wahrscheinlich"
in the German original) in connection with it. Riemann goes on to say
in his paper: ``One would of course like to have a rigorous proof of
this, but I have put aside the search for such a proof after some
fleeting vain attempts because it is not necessary for the immediate
objective of my investigation".
 The RH is undoubtedly one of the most celebrated and difficult open
problems in whole Mathematics. Its proof (or disproof) would
have very important consequences in multiplicative number theory,
especially in problems involving the distribution of primes. It would also
very likely lead to generalizations to many other zeta-functions
(Dirichlet series) having similar properties as  $\zeta(s)$.

\medskip
The RH can be put into many equivalent forms. One of the classical
is
$$
\pi(x) = {\rm li}\,x + O(\sqrt{x}\log x), \leqno{(1.5)}
$$

\noindent where $\pi(x)$ is the number of primes not exceeding
$x  \; ( \ge 2)$ and
$$
{\rm li}\,x = \int_0^x{\d t \over \log t} =
\lim_{\varepsilon \to 0+}\Bigl(\int_0^{1-\varepsilon} +
\int_{1+\varepsilon}^x \Bigr) {\d t \over {\log t}} =
\sum_{n=1}^N{(n-1)!x\over \log^nx} + O\bigl({x \over \log^{N+1}x}\bigr)
\leqno{(1.6)}
$$

\noindent
for any fixed integer $N \ge 1$. One can give a purely arithmetic
equivalent of the RH without mentioning primes. Namely we can
define recursively the M\"obius function $\mu(n)$ as
$$
\mu(1) = 1, \; \mu(n) = -\,\sum_{d\vert n, d < n}\mu(d) \qquad (n > 1).
$$

\noindent Then the RH is equivalent to the following assertion:
For any given integer $k \ge 1$ there exists an integer
$N_0 = N_0(k)$ such that, for integers $N \ge N_0$, one has
$$
{\biggl(\,\sum_{n=1}^N\mu(n)\biggr)}^{2k} \le N^{k+1}. \leqno{(1.7)}
$$
\noindent
The above definition of $\mu(n)$ is elementary and avoids primes.
A non-elementary definition of $\mu(n)$ is through the series
representation
$$
\sum_{n=1}^\infty \mu(n)n^{-s} = {1 \over \zeta(s)}\qquad (\Re{\rm e}\,s > 1),
 \leqno{(1.8)}
$$
\noindent
and an equivalent form of the RH is that (1.8) holds for $\sigma > 1/2$.
The inequality (1.7) is in fact the bound
$$
\sum_{n\le x}\mu(n)\> \ll_{\varepsilon}\> x^{{1\over 2}+\varepsilon} 
\leqno{(1.9)}
$$
\noindent
in disguise, where $\varepsilon$ corresponds to $1/(2k)$, $x$ to $N$,
and the 2$k$--th power avoids absolute values. The bound (1.9) (see [16]
and [61])
is one of the classical equivalents of the RH. The sharper bound
$$
\Big\arrowvert \sum_{n\le x}\mu(n) \Big\arrowvert < \sqrt x \qquad (x > 1)
$$
was proposed in 1897 by Mertens on the
basis of numerical evidence, and later became known in the literature
as {\it the Mertens conjecture}. It was disproved in 1985 by
A.M. Odlyzko and H.J.J. te Riele [47].

\medskip
Instead of working with the complex zeros of  $\zeta(s)$ on the
critical line it is convenient to introduce the function
$$
Z(t) = \chi^{-1/2}(\hf + it)\zeta(\hf + it), \leqno{(1.10)}
$$
\noindent
where $\chi(s)$ is given by (1.2). Since $\chi(s)\chi(1-s) = 1$ and
$\overline {\Gamma(s)} = \Gamma(\overline s)$, it follows that
$\vert Z(t)\vert = \vert\zeta({1 \over 2} + it)\vert$, $Z(t)$ is even and 
$$
\overline {Z(t)} = \chi^{-1/2}(\hf - it)\zeta(\hf - it)
= \chi^{1/2}(\hf - it)\zeta(\hf + it) = Z(t).
$$
\noindent
Hence $Z(t)$ is real if $t$ is real, and the zeros of $Z(t)$
correspond to the zeros of  $\zeta(s)$ on the critical line.
Let us denote by $0 < \gamma_1 \le \gamma_2 \le \ldots$ the
positive zeros of $Z(t)$ with multiplicities counted (all
known zeros are simple). If the RH is true, then it is known
(see [61]) that
$$
S(T) = O\Bigl({\log T \over \log \log T} \Bigr), \leqno{(1.11)}
$$

\noindent
and this seemingly small improvement over (1.4) is significant.
If (1.11) holds, then
from (1.3) one infers  that $N(T + H) - N(T) > 0$ for
$H = C/\log \log T$, suitable $C > 0$ and $T \ge T_0$.
Consequently we have the bound, on the RH,
$$
\gamma_{n+1} - \gamma_n \, \ll \, {1 \over \log \log \gamma_n} \leqno{(1.12)}
$$
\noindent for the gap between consecutive zeros
on the critical line. For some unconditional results
on $\gamma_{n+1} - \gamma_n$, see [17], [18] and [25].

\medskip
We do not know exactly what motivated Riemann to conjecture
the RH. Some mathematicians, like Felix Klein, thought that
he was inspired by a sense of general beauty and symmetry in
Mathematics. Although doubtlessly the truth of the RH would provide
such harmonious symmetry, we also know now that Riemann undertook
rather extensive numerical calculations concerning $\zeta(s)$ and
its zeros. C.L. Siegel [57] studied Riemann's unpublished notes,
kept in the G\"ottingen library. It turned out that Riemann had
computed several zeros of the zeta-function and had a deep
understanding of its analytic behaviour. Siegel provided rigorous
proof of a formula that had its genesis in Riemann's work.
It came to be known later as {\it the Riemann--Siegel formula}
(see [16], [57] and [61]) and, in a weakened form, it says that
$$
Z(t) = 2\,\sum_{n\le(t/2\pi)^{1/2}}n^{-1/2}\cos\Bigl(
t\,\log {\sqrt {t/2\pi}\over n} - {t\over 2} - {\pi \over 8}\Bigr)
+ O(t^{-1/4}), \leqno{(1.13)}
$$

\noindent
where the O-term in (1.13) is actually best possible, 
namely it is $\Omega_\pm(t^{-1/4})$. 
 As usual
 $f(x) = \Omega_\pm(g(x)) $ (for $g(x) > 0$ when
$x \ge x_0$) means that we have $f(x) = \Omega_+(x)$ and
$f(x) = \Omega_-(x)$, namely that both
$$
   \limsup_{x\to \infty}{f(x) \over g(x)} > 0, \quad \liminf_{x \to
   \infty} {f(x) \over g(x)} < 0
$$
are true.  The Riemann--Siegel
formula is an indispensable tool in the theory of $\zeta(s)$,
both for theoretical investigations and for the numerical calculations
of the zeros.

\medskip
Perhaps the most important concrete reason for believing the RH is the
impressive numerical evidence in its favour. There exists a large and rich
literature on numerical calculations involving $\zeta(s)$ and its
zeros (see [38], [44], [45], [46], [51], which contain references to
further work). This literature reflects the development of Mathematics
in general, and of Numerical analysis and Analytic number theory
in particular. Suffice to say that it is known that the first 1.5 billion
complex zeros of  $\zeta(s)$ in the upper half-plane are simple
and do have real parts equal to 1/2, as predicted by the RH.
Moreover, many large blocks of zeros of much greater height
have been thoroughly investigated, and all known zeros satisfy
the RH. However, one should be very careful in relying on numerical
evidence in Analytic number theory. A classical example for this is
the inequality $\pi(x) < {\rm li}\,x$ (see (1.5) and (1.6)), noticed
already by Gauss, which is known to be true for all $x$ for which
the functions in question have been actually computed. But the inequality
$\pi(x) < {\rm li}\,x$  is false; not only does $\pi(x) - {\rm li}\,x$
assume positive values for some arbitrarily large values of $x$,
but J.E. Littlewood [37] proved that
$$
\pi(x) = {\rm li}\,x + \Omega_\pm \left(\sqrt{x}\,{{\log \log \log x}
\over {\log x}}\right).
$$
\noindent
By extending the methods of R. Sherman Lehman [56],
H.J.J. te Riele [50] showed that $\pi(x) < {\rm li}\,x$
fails for some (unspecified) $x < 6.69\times 10^{370}$. For
values of $t$ which are this large we may hope that $Z(t)$
will also show its true asymptotic behaviour. Nevertheless, we
cannot compute by today's methods the values of $Z(t)$  for  $t$
this large, actually even $t = 10^{100}$ seems out of reach at present.
To assess why the values of $t$ where $Z(t)$ will ``really" exhibit its
true behaviour must be ``very large", it suffices to compare (1.4)
and (1.11) and note that the corresponding bounds differ by a factor
of $\log \log T$, which is a very slowly varying function.
    \medskip
Just   as there are deep reasons for believing the RH, there are
also serious grounds for doubting its truth, although the author
certainly makes no claims to possess a disproof of the RH. It is in the
folklore that  several famous mathematicians, which include
 P. Tur\'an and J.E. Littlewood, believed that the RH is not true.
 The aim of this paper is to state and analyze some of the
 arguments which cast doubt on the truth of the RH.
 In subsequent sections we
 shall deal with the Lehmer phenomenon, the Davenport-Heilbronn
zeta-function, mean value formulas on the critical line, large values
on the critical line and the distribution of zeros of a class of
convolution functions. These independent topics appear to me to
be among the most salient ones which point against the truth
of the RH. The first two of them, the Lehmer phenomenon and the
Davenport-Heilbronn zeta-function, are classical and fairly well
known. The remaining ones are rather new and are connected with the
author's research, and for these reasons the emphasis will be on
them. A sharp asymptotic formula for the convolution function 
$M_{Z,f}(t)$, related to $Z(t)$, is given in Section 8. Finally
a conditional disproof of the RH, based on the use of the
 functions $M_{Z,f}(t)$, is given at the end of the paper in
Section 9. Of course, nothing short of {\it rigorous} proof or disproof
will settle the truth of the RH.

\bigskip {\bf Acknowledgement}. I want to thank Professors M. Jutila,
K. Matsumoto, Y. Motohashi and A.M. Odlyzko for valuable remarks.

\bigskip

\centerline{\dd 2. Lehmer's phenomenon}
\bigskip
The function $Z(t)$, defined by (1.10), has a negative local maximum
$-0.52625\ldots$ at $t = 2.47575\ldots\,$. This is the only known
occurrence of a negative local maximum, while no positive local
minimum is known. {\it Lehmer's phenomenon} (named after D.H. Lehmer,
who in his works [35], [36] made significant contributions to
the subject) is the fact (see [46] for a thorough discussion)
that the graph of $Z(t)$ sometimes barely crosses the $t$--axis.
This means that the absolute value of the maximum or minimum
of $Z(t)$ between its two consecutive zeros is small. For instance,
A.M. Odlyzko found (in the version of [46] available to the author,
but Odlyzko kindly informed me that many more examples occur in the
computations that are going on now) 1976 values of $n$ such that
$\vert Z({1 \over 2}\gamma_n + {1\over 2}\gamma_{n+1})\vert < 0.0005$
in the block that he investigated. Several extreme examples are also given by
van de Lune et al. in [38]. The Lehmer phenomenon shows the delicacy
of the RH, and the possibility that a counterexample to the RH may be
found numerically. For should it happen that, for $t \ge t_0$,
$Z(t)$ attains a negative local maximum or a positive local minimum,
then the RH would be disproved. This assertion follows (see [7])
from the following
\medskip
{\bf Proposition 1.} If the RH is true, then the graph of
$Z'(t)/Z(t)$ is monotonically decreasing between the zeros
of $Z(t)$ for $t \ge t_0$.

\smallskip
Namely suppose that $Z(t)$ has a negative local maximum or a positive
local minimum between its two consecutive zeros $\gamma_n$ and
$\gamma_{n+1}$. Then $Z'(t)$ would have at least two distinct
zeros $x_1$ and $x_2$ ($x_1 < x_2$) in ($\gamma_n, \gamma_{n+1}$),
and hence so would $Z'(t)/Z(t)$. But we have
$$
{Z'(x_1) \over Z(x_1)} < {Z'(x_2)\over Z(x_2)},
$$
\noindent
which is a contradiction, since $Z'(x_1) = Z'(x_2) = 0$.

\smallskip
To prove Proposition 1 consider the function
$$
\xi(s) := {1\over 2}s(s - 1)\pi^{-s/2}\Gamma({s\over 2})\zeta(s),
$$ \noindent
so that $\xi(s)$ is an entire function of order one (see Ch. 1 of [16]),
and one has  unconditionally
$$
{\xi '(s) \over \xi(s)} = B + \sum_\rho \bigl({1\over s - \rho} +
{1 \over \rho}\bigr)    \leqno{(2.1)}
$$
with
$$
 B = \log 2 + {1\over 2}\log \pi - 1 -
{1\over 2}C_0, 
$$ 
where $\rho$ denotes complex zeros of $\zeta(s)$ and $C_0 = =\G'(1)$
is Euler's constant. By (1.2) it follows that
$$
Z(t) = \chi^{-1/2}(\hf + it)\zeta(\hf + it)
= {\pi^{-it/2}\,\Gamma({1\over 4} + {1\over 2}it)\zeta({1\over 2} + it)
\over \vert \Gamma({1\over 4} + {1\over 2}it)\vert},
$$ \noindent
so that we may write
$$
\xi(\hf + it) = -f(t)Z(t), \quad
f(t) := \hf\pi^{-1/4}(t^2 + {\txt{1\over 4}})\vert
\Gamma({\txt{1\over 4}} + \hf it)\vert.
$$
\noindent
Consequently logarithmic differentiation gives
$$
{Z'(t)\over Z(t)} = -\,{f'(t)\over f(t)} + i\,
{\xi '({1\over 2} + it)\over \xi({1\over 2} + it)}.   \leqno{(2.2)}
$$    \noindent
Assume now that the RH is true. Then by using (2.1) with
$\rho = {1\over 2} + i\gamma, s = {1\over 2} + it$ we obtain,
if $t \not= \gamma$,
$$\Bigl({i\xi'({1\over 2} + it)\over \xi({1\over 2} + it)}\Bigr)'
= -\,\sum_\gamma {1\over (t - \gamma)^2} < -C(\log \log t)^2
\quad (C > 0)
$$ \noindent
for $t \ge t_0$, since (1.12) holds. On the other hand, by using
Stirling's formula for the gamma-function and $\log \vert z\vert
= \Re{\rm e}\,\log z$, it is readily found that
$$
{\d\over \d t}\left({f'(t)\over f(t)}\right)  \;\ll\; {1\over t},
$$  \noindent
so that from (2.2) it follows that $(Z'(t)/Z(t))' < 0$
if $t \ge t_0$, which implies Proposition 1. Actually the value of
$t_0$ may be easily effectively determined and seen not to exceed
1000. Since $Z(t)$ has no positive local minimum or negative local
maximum for $3 \le t \le 1000$, it follows that the RH is false if we
find (numerically) the occurrence of a single negative local
maximum (besides the one at $t = 2.47575\ldots$) or a positive
local minimum of $Z(t)$. It seems appropriate to quote in
concluding Edwards [7], who says  that
Lehmer's phenomenon ``must give pause to even the most convinced
believer of the Riemann hypothesis".

\bigskip\bigskip

\centerline{\dd 3. The Davenport-Heilbronn zeta-function}

\bigskip
This is a zeta-function (Dirichlet series) which satisfies a functional
equation similar to the classical functional equation (1.2)
for $\zeta(s)$. It has other analogies with $\zeta(s)$, like
having infinitely many zeros on the critical line $\sigma = 1/2$,
but for this zeta-function the analogue of the RH does not
hold. This function was introduced by H. Davenport and H.
Heilbronn [6] as
$$
f(s) = 5^{-s}\bigl(\zeta(s,{1\over 5}) + \tan \theta \, \zeta(s,{2\over 5})
- \tan \theta \,\zeta(s,{3\over 5}) - \zeta(s,{4\over 5})\bigr),
\leqno{(3.1)}
$$ \noindent
where $\theta = \arctan \,(\sqrt{10 - 2\sqrt{5}} - 2)/(\sqrt{5} - 1)$ and,
for $\Re{\rm e}\,s > 1$,
$$
\zeta(s,a) = \sum_{n=0}^\infty (n + a)^{-s} \qquad (0 < a \le 1)
$$ \noindent
is the familiar {\it Hurwitz zeta-function}, defined for $\Re{\rm e}\, \, s
\le 1$
by analytic continuation. With the above choice of
$\theta$ (see [6], [32] or [61]) it can be shown that
$f(s)$ satisfies the functional equation
$$
f(s) = X(s)f(1 - s), \quad X(s) = {2\Gamma(1 - s)\cos({\pi s \over 2})
\over 5^{s-{1\over 2}}(2\pi)^{1-s}}, \leqno{(3.2)}
$$ \noindent
whose analogy with the functional equation (1.2) for $\zeta(s)$
is evident. Let $1/2 < \sigma_1 < \sigma_2 < 1$. Then it can be
shown (see Ch. 6 of [32]) that $f(s)$ has infinitely many zeros
in the strip $\sigma_1 < \sigma = \Re{\rm e}\, s < \sigma_2$, and it also has
(see Ch. 10 of [61]) an infinity of zeros in the half-plane
$\sigma > 1$, while from the product representation in (1.1) it
follows that $\zeta(s) \not= 0$ for $\sigma > 1$,
so that in the half-plane $\sigma > 1$ the behaviour of
zeros of $\zeta(s)$ and $f(s)$ is different. Actually the number
of zeros of $f(s)$ for which $\sigma > 1$ and $0 < t = \Im{\rm m}\, s \le T$
is $\gg T$, and similarly each rectangle $0 < t \le T,
1/2 < \sigma_1 < \sigma\le\sigma_2 \le 1$ contains at least
$c(\sigma_1,\sigma_2)T$
zeros of $f(s)$. R. Spira [58] found that $0.808517 + 85.699348i$
(the values are approximate)
is a zero of $f(s)$ lying in the critical strip $0 < \sigma < 1$, but
not on the critical line $\sigma = 1/2$. On the other hand,
A.A. Karatsuba [31] proved that the number of zeros
${1\over 2} + i\gamma$ of $f(s)$ for which $0 < \gamma \le T$
is at least $T(\log T)^{1/2-\varepsilon}$ for any given
$\varepsilon > 0$ and $T \ge T_0(\varepsilon)$. This bound is weaker
than A. Selberg's classical result [53] that there are
$\gg T\,\log T$ zeros ${1\over 2} + i\gamma$ of $\zeta(s)$ for which
$0 < \gamma \le T$. From the Riemann--von Mangoldt formula (1.3)
it follows that, up to the value of the $\ll$--constant, Selberg's
result on $\zeta(s)$ is best possible. There are certainly
$\ll T\,\log T$ zeros ${1\over 2} + i\gamma$ of $f(s)$
for which $0 < \gamma \le T$ and it may be that almost all of them
lie on the critical line $\sigma = 1/2$, although this has not been
proved yet. The Davenport-Heilbronn zeta-function is not the only example
of a zeta-function that exhibits the phenomena described above,
and many so-called {\it Epstein zeta-functions} also have
complex zeros off their respective critical lines (see the paper
of E. Bombieri and D. Hejhal [5] for some interesting results).

\smallskip
What is the most important difference between $\zeta(s)$ and
$f(s)$ which is accountable for the difference of distribution
of zeros of the two functions, which occurs at least
in the region $\sigma > 1$? 
It is most likely that the answer
is the lack of the Euler product for $f(s)$, similar to the one
in (1.1) for $\zeta(s)$.  But $f(s)$ can be written as a linear
combination of two $L$-functions which have Euler products
(with a common factor) and this fact plays the crucial r\^ole in
Karatsuba's proof of the lower bound result for the number of zeros
of $f(s)$. In any case one can argue that it may
likely happen that the influence of the Euler product for $\zeta (s)$
will not extend all the way to the line $\sigma = 1/2$.
In other words, the existence of zeta-functions
such as $f(s)$, which share many common properties with $\zeta(s)$,
but which have infinitely many zeros off the critical line, certainly
disfavours the RH.

\medskip
Perhaps one should at this point mention the {\it Selberg
zeta-function} ${\cal Z}(s)$ (see [55]). This is an entire
 function which enjoys several
common properties with $\zeta (s)$, like the functional equation and
the Euler product. For ${\cal Z}(s)$ the corresponding analogue of the RH is
true,
 but it should be stressed that 
${\cal Z}(s)$ is {\it not} a classical Dirichlet series. Its Euler product
$$
{\cal Z}(s) = \prod_{\{P_0\}}\prod_{k=0}^\infty(1 - N(P_0)^{-s-k}) \qquad
(\Re{\rm e}\,s > 1)
$$
\noindent
is not a product over the rational primes, but over norms of certain
conjugacy classes of groups. Also ${\cal Z}(s)$ is an entire function
 of order 2, while $(s-1)\zeta(s)$ is an entire function of order 1.
For these reasons ${\cal Z}(s)$ cannot be compared
too closely to $\zeta (s)$.

\bigskip
\centerline{\dd 4. Mean value formulas on the critical line}
\bigskip

For $k \ge 1$ a fixed integer, let us write the $2k$-th moment of $|\zt|$ as
$$
\int_0^T\vert\zeta(\hf + it)\vert^{2k}\d t =
T\,P_{k^2}(\log T) + E_k(T), \leqno{(4.1)}
$$
\noindent
where for some suitable coefficients $a_{j,k}$ one has
$$
P_{k^2}(y) = \sum_{j=0}^{k^2}a_{j,k}y^j.  \leqno{(4.2)}
$$ \noindent
An extensive literature exists on $E_k(T)$, especially on
$E_1(T) \equiv E(T)$ (see F.V. Atkinson's classical paper [2]),
 and the reader is referred to [20] for a
 comprehensive account. It is known that
 $$
 P_1(y) = y + 2C_0 - 1 - \log(2\pi),
 $$ \noindent
 and $P_4(y)$ is a quartic polynomial whose leading coefficient
 equals $1/(2\pi^2)$ (see [22] for an explicit evaluation of its
 coefficients). One hopes that
$$
E_k(T) = o(T) \qquad (T \to \infty) \leqno{(4.3)}
$$ \noindent
will hold for each fixed integer $k \ge 1$, but so far this is
known to be true only in the cases $k = 1$ and $k = 2$, when
$E_k(T)$ is a true error term in the asymptotic formula (4.1).
In fact heretofore it has not been clear how to define properly (even
on heuristic grounds) the values of $a_{j,k}$ in (4.2)  for
$k \ge 3$ (see [24] for an extensive discussion concerning the
case $k = 3$). The connection between $E_k(T)$ and the RH is indirect,
namely there is a connection with the {\it Lindel\"of hypothesis}
(LH for short). The LH is also a famous unsettled problem, and it
states that
$$
\zeta(\hf + it) \> \ll_\varepsilon\>t^\varepsilon \leqno{(4.4)}
$$ \noindent
for any given $\varepsilon > 0$ and $t \ge t_0 > 0$ (since
$\overline {\zeta({1\over 2} + it)} = \zeta({1\over 2} - it)$,
$t$ may be assumed to be positive).
It is well-known (see  [61] for a proof) that the RH implies
$$
\zeta(\hf + it) \>\ll\> \exp\bigl({A\,\log t \over \log \log t}
\bigr) \qquad (A > 0,\; t \ge t_0), \leqno{(4.5)}
$$ \noindent
so that obviouly the RH implies the LH. In the other direction it is
unknown whether the
LH (or (4.5)) implies the RH.
However, it is known that the LH has considerable influence on the
distribution of zeros of $\zeta (s)$. If $N(\sigma ,T)$ denotes the
number of zeros $\rho = \beta  + i\gamma $ of $\zeta (s)$
for which $\sigma  \le \beta $ and $\vert\gamma\vert \le T$,
then it is known (see Ch. 11 of [16]), that the LH  implies  that
$N(\sigma,T) \ll T^{2-2\sigma+\varepsilon}$ for $1/2 \le \sigma  \le 1$
(this is a form of {\it the density hypothesis}) and
$N({3\over 4} + \delta, T) \ll T^\varepsilon$, where $\varepsilon =
\varepsilon(\delta)$ may be arbitrarily small for any $0 < \delta <
{1\over 4}$.

\smallskip

 The best unconditional bound for the order of $\zeta (s)$ on the
 critical line, known at
the time of the writing of this text is
$$
\zeta(\hf + it) \>\ll_\varepsilon \>t^{c+\varepsilon} \leqno{(4.6)}
$$ \noindent
with $c = 89/570 = 0.15614\ldots\;$. This is due to M.N. Huxley [13],
and represents the last in a long series of improvements over
the past 80 years. The result is obtained by intricate estimates of
exponential sums of the type $\sum_{N<n\le 2N}n^{it} \quad (N \ll \sqrt t)$,
and the value $c = 0.15$ appears to be the limit of the method.

\medskip
Estimates for $E_k(T)$ in (4.1) (both pointwise and in the mean sense)
have many applications. From the knowledge about the order of $E_k(T)$
one can deduce a bound for $\zeta({1\over 2} + iT)$ via the estimate
$$
\zeta(\hf + iT) \>\ll\> (\log T)^{(k^2+1)/(2k)} +
\bigl(\,\log T \max_{t\in [T-1,T+1]} \vert E_k(t)\vert\,\bigr)^{1/(2k)},
\leqno{(4.7)}
$$ \noindent
which is Lemma 4.2 of [20].  Thus the best known upper bound
$$E(T) \equiv E_1(T) \ll T^{72/227}(\log T)^{679/227}
\leqno{(4.8)}
$$

\noindent
of M.N. Huxley [14] yields (4.6) with $c = 36/227 = 0.15859\ldots\;$.
Similarly the sharpest known bound
$$
E_2(T) \ll T^{2/3}\log ^CT \qquad (C > 0) \leqno{(4.9)}
$$ \noindent
of Y. Motohashi and the author (see [20], [26], [28]) yields
(4.6) with the classical value $c = 1/6$ of Hardy and Littlewood.
Since the difficulties in evaluating the left-hand side of (4.1)
greatly increase as $k$ increases, it is reasonable to expect that
the best estimate for $\zeta({1\over 2} + iT)$ that one can get
from (4.7) will be when $k = 1$.

\medskip
The LH is equivalent to the bound
$$\int_0^T\vert\zeta(\hf + it)\vert^{2k}\d t \ll_{k,\varepsilon}
T^{1+\varepsilon}   \leqno{(4.10)}
$$ \noindent
for any $k \ge 1$ and any $\varepsilon > 0$, which in turn is the same as
$$
E_k(T) \ll_{k,\varepsilon} T^{1+\varepsilon}. \leqno{(4.11)}
$$ \noindent
The enormous difficulty in settling the truth of the LH, and so
{\it a fortiori} of the RH, is best reflected in the relatively modest
upper bounds for the integrals in (4.10) (see Ch. 8 of [16] for
sharpest known results). On the other hand, we have $\Omega$-results
in the case $k = 1,2$, which show that $E_1(T)$ and $E_2(T)$
cannot be always small. Thus J.L. Hafner and the author [11], [12]
proved that
$$
E_1(T) = \Omega_+\Bigl((T\,\log T)^{1\over 4}(\log \log T)^{3+\log4\over 4}
{\rm e}^{-C\sqrt{\log \log \log T}}\Bigr) \leqno{(4.12)}
$$ \noindent
and
$$
E_1(T) = \Omega_-\biggl(T^{1\over 4}\exp\Bigl({D(\log \log T)^{1\over 4}
\over (\log \log \log T)^{3\over 4}}\Bigr)\biggr) \leqno{(4.13)}
$$
\noindent
for some absolute constants $C,D >0$. Moreover the author [19]
proved that there exist constants $A,B > 0$ such that, for $T \ge T_0$,
every interval $[T, T + B\sqrt{T}\,]$ contains points $t_1, t_2$ for which
$$
E_1(t_1) > At_1^{1/4}, \quad E_1(t_2) < -At_2^{1/4}.
$$ \noindent
Numerical investigations concerning $E_1(T)$ were carried out
by H.J.J. te Riele and the author [29].

\medskip
The $\Omega$--result
$$ E_2(T) \;=\; \Omega(\sqrt{T}\,) \leqno{(4.14)}
$$
\noindent (meaning $\lim_{T\to\infty}E_2(T)T^{-1/2} \not= 0$)
was proved by Y. Motohashi and the author (see [26], [28] and Ch. 5
of [20]). The method of proof involved differences of values
of the functions $E_2(T)$, so that (4.14) was the limit of the
method. The basis of this, as well of other recent investigations
involving $E_2(T)$, is Y. Motohashi's fundamental explicit formula
for
$$
(\Delta\sqrt{\pi})^{-1}\int_{-\infty}^\infty\vert\zeta(\hf +
it + iT)\vert^4\,{\rm e}^{-(t/\Delta)^2}\d t\qquad(\Delta > 0),
\leqno{(4.15)}
$$ \noindent
obtained by deep methods involving spectral theory of the
non-Euclidean Laplacian (see [40], [41], [43], [63] and
Ch. 5 of [16]). On p. 310 of [20] it was pointed out that a
stronger result than (4.14), namely
$$
\limsup_{T\to\infty}\vert E_2(T)\vert T^{-1/2} \;=\; +\infty
$$ \noindent
follows if certain quantities connected with the discrete spectrum
of the non-Euclidean Laplacian are linearly independent over the
integers. Y. Motohashi [42] recently unconditionally improved (4.14)
by showing that
$$ E_2(T) = \Omega_\pm(\sqrt{T}) \leqno{(4.16)}
$$

\noindent
holds. Namely he proved that the function
$$
Z_2(\xi) := \int_1^\infty\vert\zeta(\hf + it)\vert^4\,t^{-\xi}\d t,
$$
\noindent
defined initially as a function of the complex variable
$\xi$ for $\Re{\rm e}\, \xi > 1$, is meromorphic over the whole
complex plane. In the half-plane $\Re{\rm e}\, \xi > 0$ it has a pole
of order five at $\xi = 1$, infinitely many simple poles of
the form ${1\over 2} \pm \kappa i$, while the remaining poles
for $\Re{\rm e}\, \xi > 0$ are of the form $\rho/2, \zeta(\rho) = 0$.
Here $\kappa^2 + {1\over 4}$ is in the discrete spectrum of the
non-Euclidean Laplacian with respect to the full modular group.
By using (4.1) and integration by parts it follows that
$$
Z_2(\xi) = C + \xi\,\int_1^\infty P_4(\log t)t^{-\xi}\d t
+ \xi\,\int_1^\infty E_2(t)t^{-\xi-1}\d t \leqno{(4.17)}
$$ \noindent
with a suitable constant $C$, where the integrals are certainly
absolutely convergent for $\Re{\rm e}\, \xi > 1$ (actually the 
 second for $\Re{\rm e}\, \xi > 1/2$ in view of (4.20)). Now
(4.16) is an immediate consequence of (4.17) and the following
version of a classical result of E. Landau (see [1] for a proof).

\medskip {\bf Proposition 2}. Let $g(x)$ be a continuous function
such that
$$
G(\xi) := \int_1^\infty g(x)x^{-\xi-1}\d x
$$ \noindent
converges absolutely for some $\xi $. Let us suppose that $G(\xi)$
admits analytic continuation to a domain including the half-line
$[\sigma,\infty)$, while it has a simple pole at $\xi = \sigma +
i\delta \;(\delta \not= 0)$, with residue $\gamma$. Then
$$
\limsup_{x\to\infty}g(x)x^{-\sigma} \ge \vert \gamma\vert,
\quad \liminf_{x\to\infty}g(x)x^{-\sigma} \le -\vert \gamma\vert.
$$

It should be pointed out that (4.14) shows that the well-known
analogy between $E_1(T)$ and $\Delta_2(x)$ ( = $\Delta(x)$, the
error term in the formula for $\sum_{n\le x}d(n)$), which is 
discussed e.g., in Ch. 15 of [16], cannot be extended to general
$E_k(T)$ and  $\Delta_k(x)$. The latter function denotes the
error term in the asymptotic formula for $\sum_{n\le x}d_k(n)$,
where $d_k(n)$ is the general divisor function generated by
$\zeta^k(s)$. The LH is equivalent to either $\alpha_k \le 1/2$
($k \ge 2$) or $\beta_k = (k - 1)/(2k)$ ($k \ge 2$), where
$\alpha_k$ and $\beta_k$ are the infima of the numbers $a_k$ and
$b_k$ for which
$$
\Delta_k(x) \ll x^{a_k}, \quad \int_1^x\Delta_k^2(y)\d y \ll x^{b_k}
$$ \noindent
hold, respectively. We know that $\beta_k = (k - 1)/(2k)$ for $k = 2,3,4$,
and it is generally conjectured that $\alpha_k = \beta_k = (k - 1)/(2k)$
for any $k$. At first I thought that, analogously to the conjecture for
$\alpha_k$ and $\beta_k$, the upper bound for general $E_k(T)$ should
be of such a form as to yield the LH when $k \to \infty$, but in view
of (4.14) I am certain that this cannot be the case.

\bigskip
It may be asked then
how do the $\Omega$-results for $E_1(T)$ and $E_2(T)$ affect the LH,
and thus indirectly the RH? A reasonable conjecture is that these
$\Omega$-results lie fairly close to the truth, in other words
that$$
E_k(T) = O_{k,\varepsilon}(T^{{k\over 4}\,+\,\varepsilon}) \leqno{(4.18)}
$$ \noindent
holds for $k = 1,2$. This view is suggested by estimates in the mean
for the functions in question. Namely the author [15] proved that
$$
\int_1^T\vert E_1(t)\vert^A\d t \ll_\varepsilon  T^{1+{A\over 4}+\varepsilon}
\qquad (0 \le A \le {35\over 4}), \leqno{(4.19)}
$$ \noindent
and the range for $A$ for which (4.19) holds
can be slightly increased by using the best
known estimate (4.6) in the course of the proof. Also Y. Motohashi
and the author [27], [28] proved that
$$
\int_0^TE_2(t)\d t \ll T^{3/2},\quad \int_0^TE_2^2(t)\d t \ll
T^2\log^CT  \quad (C > 0). \leqno{(4.20)}
$$ \noindent
The bounds (4.19) and (4.20) show indeed that, in the mean sense,
the bound (4.18) does hold when $k = 1,2$. Curiously enough,
it does not seem  possible to show that the RH implies
(4.18) for $k \le 3$. If (4.18) holds for
any $k$, then in view of (4.7) we would obtain (4.6) with the
hitherto sharpest bound $c \le 1/8$, or equivalently
$\mu(1/2) \le 1/8$, where for any real $\sigma$ one defines
$$
\mu(\sigma) = \limsup_{t\to\infty}{\log\,\vert\zeta(\sigma + it)
\vert\over\log t},
$$ \noindent
and it will be clear from the context that no confusion can arise
with the M\"obius function.   What can one expect about the order
of magnitude of $E_k(T)$ for $k \ge 3$? It was already mentioned
that the structure of $E_k(T)$ becomes increasingly complex as
$k$ increases. Thus we should not expect a smaller exponent
than $k/4$ in (4.18)
for $k \ge 3$, as it would by (4.7) yield a result of the type
$\mu(1/2) < 1/8$, which in view of the $\Omega$-results is not
obtainable from (4.18) when $k = 1,2$. Hence by analogy with the cases $k = 1,2$
one would be led to conjecture that

$$
E_k(T) = \Omega(T^{k/4}) \leqno{(4.21)}
$$ \smallskip\noindent
holds for any fixed $k \ge 1$. But already for $k = 5$ \enskip (4.21)
yields, in view of (4.1),
$$
\int_0^T\vert\zeta(\hf + it)\vert^{10}\d t = \Omega_+(T^{5/4}),
\leqno{(4.22)}
$$ \noindent
which contradicts (4.10), thereby disproving {\it both the LH and the RH}.
It would be of great interest to obtain more detailed information on
$E_k(T)$ in the cases when $k = 3$ and especially when $k = 4$, as
the latter probably represents a turning point in the asymptotic
behaviour of mean values of $\vert\zeta({1\over 2} + it)\vert$.
Namely the above phenomenon strongly suggests that either the LH
fails, or the shape of the asymptotic formula for the left-hand side
of (4.1) changes (in a yet completely unknown way) when $k = 4$. In
[24] the author proved that $E_3(T) \ll_\varepsilon T^{1+\varepsilon}$
conditionally,
that is, provided that a certain conjecture involving the ternary
additive divisor problem holds. Y. Motohashi ( [40] p. 339, and [42])
 proposes,
on heuristic grounds based on analogy with explicit formulas known in
the cases $k = 1,2$, a formula for the analogue of (4.15) for the
sixth moment, and also conjectures (4.21) for $k = 3$. Concerning the
eighth moment, it should be mentioned that N.V. Kuznetsov [33]
had an interesting approach based on applications of spectral theory,
but unfortunately his proof of
$$
\int_0^T\vert\zeta(\hf+ it)\vert^8\d t \ll T\log^CT
\leqno{(4.23)}$$
\noindent
had several gaps (see the author's review in Zbl. 745.11040 and
the Addendum of Y. Motohashi [40]), so that (4.23) is still a
conjecture. If (4.23) is true, then one must have $C \ge 16$ in (4.23),
since by a result of K. Ramachandra (see [16] and [48]) one has,
for any rational number $k \ge 0$,
$$
\int_0^T\vert\zeta(\hf + it)\vert^{2k}\d t \gg_k T(\log T)^{k^2}.
$$

\medskip
The LH (see [61]) is equivalent to the statement that $\mu(\sigma)
= 1/2 - \sigma$ for $\sigma < 1/2$, and $\mu(\sigma) = 0$ for
$\sigma \ge 1/2$. If the LH is not true, what would then the graph
of $\mu(\sigma)$ look like? If the LH fails, it is most likely that $\mu(1/2)
= 1/8$ is true. Since $\mu(\sigma)$ is (unconditionally) a non-increasing,
convex function of $\sigma$,
$$
\mu(\sigma) = {1\over 2} - \sigma \quad (\sigma \le 0), \quad
\mu(\sigma) = 0 \quad(\sigma \ge 1),
$$     \noindent
and by the functional equation one has
$\mu(\sigma) = {1\over 2} - \sigma + \mu (1 - \sigma) \; (0 <
\sigma  < 1),$ perhaps one would have
$$
\mu(\sigma) = \cases{{1\over 2} - \sigma &$\sigma \le {1\over 4},$\cr  &\cr
{3\over 8} - {\sigma\over 2} &${1\over 4} \le \sigma < {3\over 4},$\cr  &\cr
0 &$\sigma \ge {3\over 4},$\cr} \leqno{(4.24)}
$$ \noindent or the slightly weaker
$$
\mu(\sigma) = \cases{{1\over 2} - \sigma &$\sigma \le 0,$\cr  &\cr
{2 - 3\sigma\over 4} &$0 < \sigma < {1\over 2},$\cr  &\cr
{(1 - \sigma)\over 4} &${1\over 2} \le \sigma \le 1,$\cr   &\cr
0 &$\sigma > 1.$\cr} \leqno{(4.25)}
$$ \noindent
A third candidate is
$$
\mu(\sigma) = \cases{{1\over 2} - \sigma &$\sigma \le 0,$\cr &\cr
{1\over 2}(1 - \sigma)^2 &$0 < \sigma < 1,$\cr &\cr
0 &$\sigma \ge 1,$\cr} \leqno{(4.26)}
$$ \noindent
which is a quadratic function of $\sigma$ in the critical strip.
 Note that (4.26)
sharpens (4.25) for $0 < \sigma < 1$, except when $\sigma = 1/2$, when
(4.24)--(4.26) all yield $\mu (\sigma) = 1/8$.
So far no exact value of $\mu(\sigma)$ is known when $\sigma$
lies in the critical strip $0 < \sigma < 1$.

\bigskip   \medskip
\centerline{\dd 5. Large values on the critical line}
\bigskip
One thing that has constantly made the author skeptical about the
truth of the RH is: How to draw the graph of $Z(t)$ when $t$ is
large? By this the following is meant. R. Balasubramanian
and K. Ramachandra (see [3], [4], [48], [49]) proved
unconditionally that
$$
\max_{T\le t\le T+H}\vert\zeta(\hf + it)\vert
 > \exp\Bigl({3\over 4}\Bigl({\log H\over \log\log H}\Bigr)^{1/2}
 \Bigr) \leqno{(5.1)}
 $$ \noindent
for $T \ge T_0$ and $\log\log T \ll H \le T$, and probably
on the RH this can be further improved (but no results seem to
exist yet). Anyway (5.1) shows that $\vert Z(t)\vert$ assumes
large values relatively often. On the other hand, on the RH
one expects that the bound in (1.11) can be also further reduced,
very likely (see [46]) to
$$
S(T) \ll_\varepsilon (\log T)^{{1\over 2} + \varepsilon}.
\leqno{(5.2)}
$$ \noindent
Namely, on the RH, H.L. Montgomery [39] proved that
$$
S(T) = \Omega_\pm\Bigl(\Bigl( {\log T \over \log \log T }\Bigr)^{1/2}\Bigr),
$$ \noindent
which is in accord with (5.2). K.-M. Tsang [59], improving a
classical result of A. Selberg [54], has shown that one has
unconditionally
$$
S(T) = \Omega_\pm\Bigl(\Bigl( {\log T \over \log \log T }\Bigr)^{1/3}\Bigr).
$$ \noindent
Also K.-M. Tsang [60] proved that (unconditionally; $\pm$ means that
the result holds both with the $+$ and the $-$ sign)
$$
\bigl(\sup_{T\le t \le 2T}\log\vert\zeta(\hf + it)\vert\bigr)
\bigl(\sup_{T\le t \le 2T}\pm S(t)\bigr) \gg {\log T\over \log \log T},
$$ \noindent
which shows that either $\vert\zeta({1\over 2} + it)\vert$ or
$\vert S(t)\vert$ must assume large values in $[T, 2T]$. It may be
pointed out that the calculations relating to the values of $S(T)$
(see e.g., [45], [46]) show that all known values of $S(T)$ are
relatively small. In other words they are not anywhere near the values predicted
by the above $\Omega$--results, which is one more reason that supports
the view that the values for which $\zeta (s)$ will exhibit its
true asymptotic behaviour must be {\it really very large}.

\medskip If on the RH (5.2) is true, then clearly (1.12) can be
improved to
$$
\gamma_{n+1} - \gamma_n \ll_\varepsilon (\log \gamma_n)^{\varepsilon - 1/2}.
\leqno{(5.3)}
$$ \noindent
This means that, as $n \to \infty$, the gap between the consecutive
zeros of $Z(t)$ tends to zero not so slowly. Now take $H = T$ in
(5.1), and let $t_0$ be the point in $[T, 2T]$ where the maximum
in (5.1) is attained. This point falls into an interval of length
$\ll (\log T)^{\varepsilon - 1/2}$ between two consecutive zeros,
so that in the vicinity of $t_0$ the function $Z(t)$ must have very
large oscillations, which will be carried over to $Z'(t), Z''(t),\ldots$
etc. For example, for $T = 10^{5000}$ we shall have
$$
\vert Z(t_0)\vert > 2.68\times 10^{11}, \leqno{(5.4)}
$$ \medskip \noindent
while $(\log T)^{-1/2} = 0.00932\ldots\,$ , which shows how large
the oscillations of $Z(t)$ near $t_0$ will be. Moreover, M. Jutila
[30] unconditionally proved the following

\smallskip {\bf Proposition 3}. There exist positive constants
$a_1, a_2$ and $a_3$ such that, for $T \ge 10$, we have
$$
\exp(a_1(\log \log T)^{1/2}) \le \vert Z(t)\vert
\le \exp(a_2(\log \log T)^{1/2})
$$
\noindent
in a subset of measure at least $a_3T$ of the interval $[0, T]$.

\medskip
For $T = 10^{5000}$ one has $e^{(\log \log T)^{1/2}} =
21.28446\ldots\,$ ,
and Proposition 3 shows that relatively large values of $\vert Z(t)\vert$
are plentiful, and in the vicinity of the respective $t$'s again
$Z(t)$ (and its derivatives) must oscillate a lot. The RH and (5.3)
imply that,
as $t\to\infty$, the graph of $Z(t)$ will consist of tightly packed
spikes, which  will be more and more condensed as $t$
increases, with larger and large oscillations. This I find
hardly conceivable. Of course,
it could happen that the RH is true and that (5.3) is not.

\bigskip
\centerline{\dd 6. A class of convolution functions}

\bigskip
It does not appear easy to put the discussion of Section 5 into a
quantitative form. We shall follow now the method developed by
the author in [21] and [23] and try to make a self-contained presentation,
resulting in the proof of  Theorem 1 (Sec. 8) 
and Theorem 2 (Sec. 9). The basic idea is to connect the order of
$Z(t)$ with the distribution of its zeros and the order of its
derivatives (see (7.5)).
However it  turned out that if one
works directly with $Z(t)$, then one encounters several difficulties. 
One is that we do not know yet whether the zeros of $Z(t)$ are all distinct
(simple), even on the RH (which implies by (1.11) only the fairly
weak bound that the multiplicities of zeros up to height $T$ are
$\ll \log T/\log \log T$). This difficulty is technical, and we may
bypass it by using a suitable form of divided differences from
Numerical analysis, as will be shown a little later in Section 7.
A.A. Lavrik [34] proved the useful result
that, uniformly for $0 \le k \le {1\over 2}\log t$, one has
$$
Z^{(k)}(t) = 2\sum_{n\le(t/2\pi)^{1/2}}n^{-1/2}
{\bigl(\log\,{(t/2\pi)^{1/2}\over n}\bigr)}^k
\cos\biggl(t\log\,{(t/2\pi)^{1/2}\over n} - {t\over 2} - {\pi\over 8}
+ {\pi k\over 2}\biggr) +
O\left(t^{-1/4}{({\txt{3\over 2}}\log t)}^{k+1}\right).\leqno{(6.1)}
$$         \noindent
The range for which (6.1) holds is large, but it is difficult to obtain
good uniform bounds for $Z^{(k)}(t)$ from (6.1). To overcome this
obstacle the author introduced in [21] the class of convolution functions
$$
M_{Z,f}(t): = \int_{-\infty}^\infty Z(t - x)f_G(x)\d x =
 \int_{-\infty}^\infty Z(t + x)f({x\over G})\d x, \leqno{(6.2)}
$$ \noindent
where $G > 0, f_G(x) = f(x/G)$, and $f(x) \; ( \ge 0)$ is an even
function belonging to the class of smooth ($C^\infty $) functions
$f(x)$ called $S_\alpha^\beta$ by Gel'fand and Shilov [9]. The
functions $f(x)$ satisfy for any real $x$ the inequalities
$$
\vert x^kf^{(q)}(x)\vert  \le CA^kB^qk^{k\alpha}q^{q\beta}
\qquad (k,q = 0,1,2,\ldots) \leqno{(6.3)}
$$ \noindent
with suitable constants $A,B,C > 0$ depending on $f$ alone. For
$\alpha = 0$ it follows that $f(x)$ is of bounded support,
namely it vanishes for $\vert x\vert \ge A$. For $\alpha > 0$
the condition (6.3) is equivalent (see [9]) to the condition
$$
\vert f^{(q)}(x)\vert \le CB^qq^{q\beta}\exp(-a\vert x\vert^{1/\alpha})
\qquad (a = \alpha/({\rm e}A^{1/\alpha})) \leqno{(6.4)}
$$ \noindent
for all $x$ and $q \ge 0$. We shall denote by $E_\alpha^\beta$
the subclass of $S_\alpha^\beta$ with $\alpha > 0$ consisting of
 even functions $f(x)$ such that $f(x)$ is not the
zero-function. It is shown in [9] that $S_\alpha^\beta$
is non-empty if $\beta \ge 0$ and $\alpha + \beta \ge 1$. If \smallskip
\noindent
these conditions hold then $E_\alpha^\beta$ is also non-empty, since
$f(-x) \in S_\alpha^\beta$ if  $f(x) \in S_\alpha^\beta$, and
$f(x) + f(-x)$ is always even.

\medskip
One of the main properties of the convolution function $M_{Z,f}(t)$,
which follows by $k$-fold integration by parts from (6.2), is that for any
integer $k \ge 0$
$$
M_{Z,f}^{(k)}(t) = M_{Z^{(k)},f}(t) =
\int_{-\infty}^\infty Z^{(k)}(t + x)f({x\over G})\d x
= {\Bigl({-1 \over G}\Bigr)}^k\,\int_{-\infty}^\infty Z(t + x)f^
{(k)}({x\over G})\d x. \leqno{(6.5)}
$$ \noindent
This relation shows that the order of $M^{(k)}$ depends only on
the orders of $Z$ and $f^{(k)}$, and the latter is by (6.4) of
exponential decay, which is very useful in dealing with convergence
problems etc.  The  salient
point of our approach is that the difficulties inherent in the
distribution of zeros of $Z(t)$ are transposed to the
distribution of zeros of $M_{Z,f}(t)$, and for the latter function
(6.5) provides good uniform control of its derivatives.

\smallskip
Several analogies between $Z(t)$  and $M_{Z,f}(t)$
are established in [21], especially in connection with mean values
and the distribution of their respective zeros. We shall retain here
the notation introduced in [21], so that $N_M(T)$ denotes the number
of zeros of $M_{Z,f}(t)$ in $(0,T]$, with multiplicities counted.
If $f(x) \in E_\alpha^\beta$, $f(x) \ge 0$ and $G = \delta/\log(T/(2\pi))$
with suitable $\delta > 0$, then Theorem 4 of [21] says that
$$
N_M(T + V) - N_M(T - V) \gg {V\over \log T}, \quad V = T^{c+\varepsilon},
\quad c = 0.329021\ldots, \leqno{(6.6)}
$$
for any given $\varepsilon > 0$. The nonnegativity of $f(x)$ was needed
in the proof of this result. For the function $Z(t)$ the analogous result
is that
$$
N_0(T + V) - N_0(T - V) \gg V\log T, \quad  V = T^{c+\varepsilon},
\quad c = 0.329021\ldots, \leqno{(6.7)}
$$
where as usual $N_0(T)$ denotes the number of zeros of $Z(t)$ (or of
$\zeta({1\over 2} + it)$)  in $(0,T]$, with multiplicities counted.
Thus the fundamental problem in the theory of $\zeta(s)$ is to
estimate $N(T) - N_0(T)$, and
the RH may be reformulated as $N(T) = N_0(T)$ for $T > 0$. The
bound (6.7) was proved by A.A. Karatsuba (see [32] for a detailed
account). As explained in [21], the bound (6.6) probably falls short
(by a factor of $\log^2T$) from the expected (true) order of magnitude
for the number of zeros of $M_{Z,f}(t)$ in $[T - V,T +V]$. This is due
to the method of proof of (6.6), which is not as strong as the classical method
of A. Selberg [54] (see also Ch. 10 of [61]). The function $N_M(T)$ seems
much more difficult to handle than $N_0(T)$ or $N(T)$. The latter can
be conveniently expressed (see [16] or [61]) by means of a complex
integral from which one infers then (1.3) with the bound (1.4). I was unable
to find an analogue of the integral representation for $N_M(T)$. Note that
the bound on the right-hand side of (6.7) is actually of the best possible
order of magnitude.

In the sequel we shall need the following technical result, which we state
as

\smallskip {\bf Lemma 1}. {\it If $L = (\log T)^{{1 \over 2}+\varepsilon},
 P = \sqrt{T\over 2\pi},
0 < G < 1, L \ll V \le T^{1\over 3}, f(x) \in E_\alpha^\beta$, then}
$$
\int\displaylimits_{T-VL}^{T+VL}\vert M_{Z,f}(t)\vert
{\rm e}^{-(T-t)^2V^{-2}}\d t
\ge GV\{\vert \widehat{f}({G \over 2\pi}\log P)\vert +
 O(T^{-1/4} + V^2T^{-3/4}L^2)\}.\leqno{(6.8)}
$$

{\bf Proof}. In (6.8) $\widehat{f}
(x)$ denotes the Fourier transform of $f(x)$,
namely
$$
\widehat{f}(x) = \int_{-\infty}^\infty f(u){\rm e}^{2\pi ixu}\d u =
\int_{-\infty}^\infty f(u)\cos(2\pi xu)\d u + i\int_{-\infty}^\infty
f(u)\sin(2\pi xu)\d u
= \int_{-\infty}^\infty f(u)\cos(2\pi xu)\d u
$$
\noindent
since $f(x)$ is even. From the Riemann--Siegel formula (1.13) we have,
if $\vert T - t\vert \le VL$ and $\vert x\vert \le \log^CT \; (C > 0)$,
$$
Z(t + x) = 2\sum_{n\le P}n^{-1/2}\cos\Bigl((t+x)\log{((t+x)/(2\pi))^{1/2}
\over n} - {t + x \over 2} - {\pi\over 8}\Bigr) + O(T^{-1/4}).
$$  \noindent
Simplifying the argument of the cosine by Taylor's formula it follows that
$$
Z(t + x) = 2\sum_{n\le P}n^{-1/2}\cos\left((t+x)\log {P\over n} - {T\over 2}
- {\pi\over 8}\right) + O(T^{-1/4} + V^2T^{-3/4}L^2). \leqno{(6.9)}
$$ \noindent
Hence from (6.2) and (6.9)  we have, since
$$
\cos(\alpha + \beta) = \cos\alpha\cos\beta - \sin\alpha\sin\beta
$$ \noindent
and $f(x)$ is even,
$$
M_{Z,f}(t) = 2G\sum_{n\le P}n^{-1/2}\cos(t\,\log{P\over n} - {T\over 2}
-{\pi\over 8})
\int_{-\infty}^{\infty} f(x)\cos(Gx\,\log{P\over n})\d x +
O(GT^{-1/4} + GV^2T^{-3/4}L^2)
$$

$$
= 2G\sum_{n\le P}n^{-1/2}\cos(t\,\log{P\over n} - {T\over 2} -{\pi\over 8})
\widehat{f}({G\over 2\pi}\log{P\over n}) +  O(GT^{-1/4} + GV^2T^{-3/4}L^2).
\leqno{(6.10)}
$$ \noindent
Therefore we obtain from (6.10)
$$
\int\displaylimits_{T-VL}^{T+VL}\vert M_{Z,f}(t)\vert
{\rm e}^{-(T-t)^2V^{-2}}\d t
\;\ge\; GI + O(GVT^{-1/4} + GV^3T^{-3/4}L^2),
$$ \noindent
say, where
$$
I: = \int\displaylimits_{T-VL}^{T+VL}\Big\arrowvert\sum_{n\le
P}n^{-1/2}\widehat{f}({G\over 2\pi}
\log{P\over n})\Bigl(\exp(it\log{P\over n} - {iT\over 2} - {i\pi\over 8}) +
\exp(-it\log{P\over n} + {iT\over 2} + {i\pi\over 8})\Bigr)\Big\arrowvert
{\rm e}^{-(T-t)^2V^{-2}}\d t.
$$ \noindent
By using the fact that 
$\vert \exp\bigl(it\log P - {iT\over 2} - {i\pi\over 8}\bigr)\vert = 1$
and the classical integral
$$
\int_{-\infty}^{\infty}\exp(Ax - Bx^2)\d x = \sqrt{\pi \over
B}\,\exp\bigl({A^2\over 4B}\bigr) \qquad (\Re{\rm e}\, B > 0)
$$ \noindent 
we shall obtain
$$
I \;\ge\; \vert I_1 + I_2\vert,
$$\noindent
where
$$ \eqalign{
I_1  &= \int\displaylimits_{T-VL}^{T+VL}
\sum_{n\le P}n^{-1/2}\widehat{f}({G\over2\pi}\log{P\over n})
\exp(-it\log n - (T-t)^2V^{-2})\d t\cr&
= \sum_{n\le P}n^{-1/2}\widehat{f}({G\over 2\pi}
\log{P\over n})\exp(-iT\log n)
\int_{-VL}^{VL}\exp(-ix\log n - x^2V^{-2})\d x\cr&
= \sum_{n\le P}n^{-1/2}\widehat{f}({G\over 2\pi}\log{P\over n})\Big\{\sqrt{\pi}V
\exp(-iT\log n - {1\over 4}V^2\log^2n) +
O(\exp(-\log^{1+2\varepsilon}T))\Big\}\cr&
=  \sqrt{\pi}V\widehat{f}({G\over  2\pi}\log P) + O(T^{-C})   \cr}
$$
\noindent
for any fixed $C > 0$. Similarly we find that
$$
I_2 =  \int\displaylimits_{T-VL}^{T+VL}
\sum_{n\le P}n^{-1/2}\widehat{f}({G\over2\pi}\log{P\over n})
\exp\bigl(-it\log({T\over 2\pi n}) + iT + {\pi i\over 4} -
(T - t)^2V^{-2}\bigr)\d t =
$$

$$
  \sum_{n\le P}n^{-1/2}\widehat{f}({G\over 2\pi}\log{P\over n})
\exp\bigl(-it\log({T\over 2\pi n}) + iT + {\pi i\over 4}\bigr)
\Bigl(\sqrt{\pi}V\exp\{-{1\over 4}{\bigl(V\log({T\over 2\pi n})\bigr)}^2\} +
O\{\exp(-\log^{1+2\varepsilon}T)\}\Bigr)
$$
$$ = \;O(T^{-C})$$ \noindent
again for any fixed $C > 0$, since
$\log({T\over 2\pi n}) \ge \log{T\over 2\pi P} = {1\over 2}\log({T\over 2
\pi})$.
From the above estimates (6.8) follows.

\bigskip
\centerline{\dd 7. Technical preparation}

\bigskip
In this section we shall lay the groundwork for the investigation of the
distribution of zeros of $Z(t)$ via the convolution functions
$M_{Z,f}(t)$. To do this we shall first briefly outline a method based
on a generalized form of the mean value theorem from the differential
calculus. This can be conveniently obtained from the expression for the
$n$-th divided difference associated to the function $F(x)$, namely
$$
[x,x_1,x_2,\cdots,x_n]  :=
$$
$$= {F(x)\over (x-x_1)(x-x_2)\cdots(x-x_n)} +
{F(x_1)\over (x_1-x)(x_1-x_2)\cdots(x_1-x_n)} + \cdots +
{F(x_n)\over (x_n-x)(x_n-x_1)\cdots(x_n-x_{n-1})}
$$ \noindent
where $x_i \not = x_j$ if $i \not = j$, and $F(t)$ is a real-valued
function of the real variable $t$. We have the representation
$$
[x,x_1,x_2,\cdots,x_n] =   \leqno{(7.1)}
$$
$$
=  \int_0^1\int_0^{t_1}\cdots\int_0^{t_{n-1}}F^{(n)}\bigl(x_1 +
(x_2-x_1)t_1 +   \cdots + (x_n-x_{n-1})t_{n-1} +
(x-x_n)t_n\bigr)\d t_n\cdots \d t_1 =
{F^{(n)}(\xi)\over n!}
$$ \noindent
if $F(t) \in C^n[I], \xi = \xi(x,x_1,\cdots,x_n)$ and $I$ is the
smallest interval containing all the points $x,x_1,\cdots,x_n$. If
we suppose additionally that $F(x_j) = 0$ for $j = 1,\cdots,n$, then on
comparing the two expressions for $[x,x_1,x_2,\cdots,x_n]$, it follows that
$$
F(x) = (x - x_1)(x - x_2)\cdots(x - x_n){F^{(n)}(\xi)\over n!},\leqno{(7.2)}
$$ \noindent
where $ \xi = \xi(x)$ if we consider $x_1,\cdots,x_n$ as fixed and $x$ as
a variable. The  underlying idea is that, if the (distinct) zeros $x_j$
of $F(x)$ are sufficiently close to one another, then (7.2) may lead to
a contradiction if $F(x)$ is assumed to be large and one has good bounds for
its derivatives.

\medskip
To obtain the analogue of (7.2) when the points $x_j$ are not necessarily
distinct, note that if $F(z)$ is a regular function of the complex
variable $z$ in a region which contains the distinct points
$x,x_1,\cdots,x_n$, then for a suitable closed contour $\cal C$ containing
these points one obtains by the residue theorem

$$
[x,x_1,x_2,\cdots,x_n]  = {1\over 2\pi i}\int\displaylimits_{\cal C}
{F(z)\over (z - x)(z - x_1)\cdots(z - x_n)}\d z.
$$\noindent
A comparison with (7.1) yields then
$$
 {1\over 2\pi i}\int\displaylimits_{\cal C}{F(z)\over
 (z - x)(z - x_1)\cdots(z - x_n)}\d z
\leqno{(7.3)}$$
$$=  \int_0^1\int_0^{t_1}\cdots\int_0^{t_{n-1}}F^{(n)}\bigl(x_1 +
(x_2-x_1)t_1 +
\cdots + (x_n-x_{n-1})t_{n-1} + (x-x_n)t_n\bigr)\d t_n\cdots \d t_1.
$$
\noindent
Now (7.3) was derived on the assumption that the points $x,x_1,\cdots,x_n$
are distinct. But as both sides of (7.3) are regular functions of
$x,x_1,\cdots,x_n$  in some region, this assumption may be dropped by
analytic continuation. Thus let the points $x,x_1,\cdots,x_n$ coincide
with the (distinct) points $z_k$, where the multiplicity of $z_k$ is
denoted by $p_k; k = 0,1,\cdots,\nu;\, \sum_{k=0}^{\nu} p_k = n + 1$.
If we set
$$
z_0 = x,\> p_0 = 1,\> Q(z) = \prod_{k=1}^{\nu}(z - z_k)^{p_k},
$$ \noindent
then the complex integral in (7.3) may be evaluated by the residue
theorem (see Ch.1 of A.O. Gel'fond [10]). It equals
$$
{F(x)\over Q(x)} - \sum_{k=1}^\nu \sum_{m=0}^{p_k-1}
{F^{(p_k-m-1)}(z_k)\over (p_k - m - 1)!}\sum_{s=0}^m{1\over (m - s)!}\,
{d^{m-s}\over dz^{m-s}}\,\Bigl({(z - z_k)^{p_k}\over Q(z)}\Bigr)
\Big\arrowvert_{z=z_k}\cdot (x - x_k)^{-s-1}. \leqno{(7.4)}
$$
\noindent
If $x_1,x_2,\cdots,x_n$ are the zeros of $F(z)$, then
$F^{(p_k-m-1)}(z_k) = 0$ for $k = 1,\cdots,\nu$ and
$m = 0,\cdots,p_{k-1}$, since $p_k$ is the multiplicity of (the zero) $z_k$.
Hence if $F(x) \not= 0$, then on comparing (7.1), (7.3) and (7.4) one
obtains
$$
\vert F(x)\vert  \;\le\;
\prod_{k=1}^n\vert x - x_k\vert{\vert F^{(n)}(\xi)\vert
\over n!} \qquad(\xi = \xi(x)),\leqno{(7.5)}
$$\noindent
and of course (7.5) is trivial if $F(x) = 0$.

\smallskip
Now we shall apply (7.5) to $F(t) = M_{Z,f}(t), f(x) \in E_\alpha^\beta $,
with $n$ replaced by $k$, to obtain
$$
\vert  M_{Z,f}(t)\vert \le \prod_{t-H  \le \gamma\le t+H}
\vert\gamma - t\vert{\vert  M_{Z,f}^{(k)}(\tau)\vert \over k!}, \leqno{(7.6)}
$$ \noindent
where $\gamma$ denotes the zeros of $ M_{Z,f}(t)$ in $[t - H, t + H],
\tau = \tau(t,H) \in [t - H, t + H],
\vert t - T\vert \le T^{1/2+\varepsilon}$
and $k = k(t,H)$ is the number of zeros of $M_{Z,f}(t)$ in $[t - H, t + H]$.
We shall choose
$$
H = {A\log_3T\over \log_2T}\qquad
(\log_rT = \log(\log_{r-1}T),\;\log_1T \equiv \log T) \leqno{(7.7)}
$$                     \noindent
for a sufficiently large $A > 0$. One intuitively feels that,
with a suitable choice (see (8.1) and (8.2)) of $G$ and $f$, the
functions $N(T)$ and $N_M(T)$ will not differ by much. Thus  we
shall suppose that the analogues of (1.3) 
and (1.11) hold for $N_M(T)$, namely that
$$
N_M(T) = {T\over 2\pi}\log({T\over 2\pi}) -  {T\over 2\pi} + S_M(T) + O(1)
\leqno{(7.8)}
$$\noindent
with a continuous function $S_M(T)$ satisfying
$$
S_M(T) = O\Bigl({\log T\over \log\log T}\Bigr), \leqno{(7.9)}
$$\noindent
although it is hard to imagine what should be the appropriate analogue
for $S_M(T)$ of the defining relation $S(T) = {1\over \pi}
\arg \zeta({1\over 2} + iT)$ in (1.4). We also suppose that
$$
\int_T^{T+U}(S_M(t + H) - S_M(t - H))^{2m}\d t \;\ll\;
U(\log(2 + H\log T))^m \leqno{(7.10)}
$$ \noindent
holds for any fixed integer $m \ge 1, T^a < U \le T,
1/2 < a \le 1, 0 < H < 1$.
Such a result holds unconditionally (even in the form of an asymptotic
formula) if $S_M(T)$ is replaced by $S(T)$, as shown in the works of
A. Fujii [8] and K.-M. Tsang [59]. Thus it seems plausible
that (7.10) will also hold. It was already mentioned that it
is reasonable to expect that $S(T)$ and $S_M(T)$ will be close to
one another. One feels that this ``closeness" should hold also in the
mean sense, and that  instead of (7.10) one could impose
a condition which links directly $S_M(T)$ and $S(T)$, such as that
for any fixed integer $m \ge 1$ one has
$$
\int_T^{T+U}(S_M(t) - S(t))^{2m}\d t \;\ll\;
U(\log\log T)^m  \qquad (T^a < U \le T,\; \hf < a \le 1).\leqno{(7.11)}
$$ \noindent
If (7.8) holds, then
$$
k = N_M(t + H) - N_M(t - H) + O(1) \leqno{(7.12)}
$$
$$= {(t + H)\over 2\pi}\log\bigl({t + H\over 2\pi}\bigr) -
{t + H\over 2\pi} - {t - H\over 2\pi}\log\bigl({t - H\over 2\pi}\bigr)
+  {t - H\over 2\pi} + S_M(t + H) - S_M(t - H) +  O(1)
$$
$$= {H\over\pi}\log ({T\over 2\pi}) + S_M(t + H) - S_M(t - H) + O(1).
$$
\smallskip To bound from above the product in (7.6) we proceed as follows.
First we have trivially
$$
\prod_{\vert\gamma - t\vert\,\le\,1/\log_2T}\vert\gamma - t\vert \le 1.
$$ \noindent
The remaining portions of the product with $t - H \le \gamma < t - 1/\log_2T$
and $t + 1/\log_2T < \gamma \le t + H$ are treated analogously, so we shall
consider in detail only the latter. We have
$$\eqalign{&
\log\Bigl(\prod_{t+1/\log_2T<\gamma\le{t+H}}\vert\gamma - t\vert\Bigr) \cr&
= \sum_{t+1/\log_2T<\gamma\le{t+H}}\log(\gamma - t)
=  \int\displaylimits_{t+1/\log_2T+0}^{t+H}\log(u - t)\d N_M(u)
\cr& =
{1\over 2\pi}\int\displaylimits_{t+1/\log_2T+0}^{t+H}
\log(u-t)\log({u\over 2\pi})\,du +
\int\displaylimits_{t+1/\log_2T+0}^{t+H}\log(u-t)\d(S_M(u)+O(1)).   \cr}
$$  \noindent
By using integration by parts and (7.9) it follows that
$$
\int\displaylimits_{t+1/\log_2T+0}^{t+H}\log(u - t)\d(S_M(u) + O(1))  =
O\Bigl({\log T\log_3T\over \log_2T}\Bigr) -
\int\displaylimits_{t+1/\log_2T}^{t+H}{S_M(u) + O(1)\over u - t}\d u
\ll {\log T\log_3T\over \log_2T},
$$ \noindent
and we have
$$\eqalign{&
{1\over 2\pi}\int\displaylimits_{t+1/\log_2T+0}^{t+H}
\log(u - t)\log({u\over 2\pi})\d u
= {1\over 2\pi}\log({T\over 2\pi})\cdot(1 + O(T^{\varepsilon-1/2}))
\int\displaylimits_{t+1/\log_2T}^{t+H}\log(u - t)\d u\cr&
=  {1\over 2\pi}\log({T\over 2\pi})\cdot\left(H\log H - H +
O\bigl({\log_3T\over\log_2T}
\bigr)\right).   \cr}
$$\noindent
By combining the above estimates we obtain

\bigskip
{\bf Lemma 2}. {\it Suppose that} (7.8) {\it and}
(7.9) {\it hold. If $\gamma$ denotes
zeros of $M_{Z,f}(t)$, $H$ is given by} (7.7) {\it and
$\vert T - t\vert \le T^{1/2+\varepsilon}$, then}
$$
\prod_{t-H\,\le\gamma\,\le t+H}\vert\gamma - t\vert \;\le\;
\exp\Big\{{1\over\pi}\log({T\over 2\pi})\cdot\left(H\log H - H +
O\left({\log_3T\over\log_2T}\right)\right)\Big\}.\leqno{(7.13)}
$$
\bigskip

\centerline{\dd 8. The asymptotic formula for the convolution function}

\bigskip
In this section we shall prove a sharp asymptotic formula for $M_{Z,f}(t)$, 
which is given by Theorem 1. This will hold if $f(x)$
belongs to a specific subclass of functions from $E_\alpha^0$ ($\alpha
 > 1$ is fixed), and for such $M_{Z,f}(t)$ we may hope that (7.8)--(7.11)
 will hold. To construct this subclass of functions
first of all let let $\varphi(x) \ge 0$ (but $\varphi(x) \not \equiv 0$)
belong to
$E_0^\alpha$. Such a choice is possible, since it is readily checked that
$f^2(x) \in S_\alpha^\beta$
if $f(x) \in S_\alpha^\beta$,
and trivially $f^2(x) \ge 0$. Thus $\varphi(x)$
is of bounded support, so that $\varphi(x) = 0$ for $\vert x\vert \ge a$
for some $a > 0$. We normalize $\varphi(x)$ so that $\int_{-\infty}^\infty
\varphi(x)\d x = 1$, and for an arbitrary constant $b > \max(1,a)$ we put
$$
\Phi(x) := \int\displaylimits_{x-b}^{x+b}\varphi(t)\d t.
$$\noindent
Then $0 \le \Phi(x) \le 1, \Phi(x)$ is even (because $\varphi(x)$ is even)
and nonincreasing for $x \ge 0$, and
$$
\Phi(x) = \cases{0\quad &${\rm if} \quad \vert x\vert \ge b + a,$\cr
          1\quad &${\rm if} \quad \vert x\vert \le b - a.$\cr}
$$ \noindent
One can also check that $\varphi(x) \in S_0^\alpha$ implies that
$\Phi(x) \in S_0^\alpha$. Namely $\vert x^k\Phi(x)\vert \le (b + a)^k$,
and for $q \ge 1$ one uses (6.3) (with $k = 0, f^{(q)}$ replaced by
$\varphi^{(q-1)},(\alpha , \beta ) = (0,\alpha) $) to obtain
$$
\vert x^k\Phi^{(q)}(x)\vert \le  (b + a)^k\vert \Phi^{(q)}(x)\vert \le 
(b + a)^k\Bigl(\vert\varphi^{(q-1)}(x +b)\vert
+ \vert\varphi^{(q-1)}(x - b)\vert\Bigr)
$$
$$
\le (b + a)^k2CB^{q-1}(q - 1)^{(q-1)\alpha} \le {2C\over B}
(b + a)^kB^qq^{q\alpha},
$$\noindent
hence (6.3) will hold for $\Phi$ in place of $f$,
with $A = b + a$ and suitable $C$. Let
$$
f(x) := \int_{-\infty}^\infty\Phi(u){\rm e}^{-2\pi ixu}\d u =
\int_{-\infty}^\infty\Phi(u)\cos(2\pi xu)\d u.
$$\noindent
A fundamental property of the class $S_\alpha^\beta$ (see [9]) is that
$\widehat{S_\alpha^\beta} = S_\beta^\alpha$, where in general
$\widehat{U} = \{\widehat{f}(x) : f(x) \in U\}$. Thus $f(x) \in S_\alpha^0$,
$f(x)$ is even (because $\Phi(x)$ is even), and by the inverse Fourier
transform we have $\widehat{f}(x) = \Phi(x)$. The function $f(x)$ is not
necessarily nonnegative, but this property is not needed in the sequel.

Henceforth  let
$$
G = {\delta\over\log({T\over 2\pi}) }\qquad(\delta > 0).\leqno{(8.1)}
$$\noindent
In view of (1.3) it is seen that, on the RH, $G$ is of the order of the
average spacing between the zeros of $Z(t)$.
If $f(x)$ is as above, then we have

\medskip
{\bf THEOREM 1}. {\it For} $\vert t - T\vert \le VL, L = \log^{{1\over
2}+\varepsilon}T,
\log^\varepsilon T \le V \le {T^{1\over 4}\over\log T}$, 
$0 < \delta < 2\pi (b - a)$ 
{\it and any fixed} $N \ge 1$ {\it we have}
$$
M_{Z,f}(t) = G\bigl(Z(t) + O(T^{-N})\bigr). \leqno{(8.2)}
$$

{\bf Proof}. Observe that the weak error term $O(T^{-1/4})$ in (8.2)
follows from (6.10) (with $x = 0$) and (6.10) when we note that
$$
\widehat{f}\Bigl({G\over 2\pi}\log{P\over n}\Bigr) =
\widehat{f}\Bigl({\delta\over 2\pi}\Bigl({1\over 2} - {\log
n\over\log({T\over 2\pi})}
\Bigr)\Bigr) = 1
$$\noindent
since by construction  $\widehat{f}(x) = 1$ for $\vert x \vert < b - a$, and
$$
\Big\arrowvert {\delta\over 2\pi}\Bigl({1\over 2} - {\log
n\over\log({T\over 2\pi})}
\Bigr)\Big\arrowvert \le {\delta\over 4\pi} < b - a\qquad \Bigl(1 \le n \le P = 
\sqrt{{T\over 2\pi}}\Bigr).
$$ \noindent
Also the hypotheses on $t$ in the formulation of Theorem 1 can be relaxed.

\smallskip

In order to prove (8.2) it will be convenient to work with the 
real-valued function $\theta(t)$, defined by
$$
Z(t) = {\rm e}^{i\theta(t)}\zeta(\hf + it) = \chi^{-1/2}(\hf + it)
\zeta(\hf + it). \leqno{(8.3)}
$$
\noindent
From the functional equation (1.2) in the symmetric form
$$
\pi^{-s/2}\Gamma({s\over 2})\zeta(s) = 
\pi^{-(1-s)/2}\Gamma({1-s\over 2})\zeta(1-s)
$$\noindent
one obtains
$$
\chi^{-1/2}(\hf + it)   =
{\pi^{-it/2}\Gamma^{1/2}({1\over 4} + {it\over 2})\over
\Gamma^{1/2}({1\over 4} - {it\over 2})},
$$
\noindent
and consequently
$$
\theta(t) = \Im{\rm m}\,\log\Gamma({\txt{1\over 4}} + \hf it)
- \hf t\log\pi . \leqno{(8.4)}
$$ \noindent
We have the explicit representation (see Ch. 3 of [32])
$$
\theta(t) = {t\over 2}\log {t\over 2\pi} - {t\over 2} -
{\pi\over 8} + \Delta(t) \leqno{(8.5)}
$$ \noindent
with ($\psi(x) = x - [x] - \hf$)
$$
\Delta(t) := {t\over 4}\log(1 + {1\over 4t^2}) +
{1\over 4}\arctan{1\over 2t} + {t\over 2}\int_0^\infty
{\psi(u)\d u\over (u + {1\over 4})^2 + t^2}. \leqno{(8.6)}
$$ \noindent
This formula is very useful, since it allows one to evaluate
explicitly all the derivatives of $\theta(t)$. For $t\to\infty$
it is seen that $\Delta(t)$ admits an asymptotic expansion in terms
of negative powers of $t$, and from (8.4) and Stirling's formula
for the gamma-function it is found that ($B_k$ is the $k$-th Bernoulli number) 
$$
\Delta(t) \> \sim\> \sum_{n=1}^\infty {(2^{2n}-1)\vert B_{2n}\vert\over
2^{2n}(2n - 1)2nt^{2n-1}}, \leqno{(8.7)}
$$ \noindent
and the meaning of (8.7) is that, for an arbitrary integer $N \ge 1$, 
$\Delta(t)$ equals the sum of the first $N$ terms of the series in (8.7),
plus the error term which is $O_N(t^{-2N-1})$. In general we shall have,
for $k \ge 0$ and suitable constants $c_{k,n}$,
$$
\Delta^{(k)}(t) \> \sim \> \sum_{n=1}^\infty c_{k,n}t^{1-2n-k}.
\leqno{(8.8)}
$$

\smallskip
For complex $s$ not equal to the poles of the gamma-factors we have the
Riemann-Siegel formula (this is equation (56) of C.L. Siegel [57])
$$\eqalign{&
\pi^{-s/2}\Gamma({s\over 2})\zeta(s) \cr&= 
\pi^{-s/2}\Gamma({s\over 2})\int\displaylimits_{0\swarrow 1}
{{\rm e}^{i\pi x^2}x^{-s}
\over {\rm e}^{i\pi x} - {\rm e}^{-i\pi x}}\d x
+ \pi^{(s-1)/2}\Gamma({1 - s\over 2})\int\displaylimits_{0\searrow 1}
{{\rm e}^{-i\pi x^2}x^{s-1}\over  {\rm e}^{i\pi x} -
{\rm e}^{-i\pi x}}\d x.\cr} \leqno{(8.9)}
$$ \noindent
\smallskip
Here $0\swarrow 1$ (resp. $0\searrow 1$) denotes a straight line
which starts from infinity in the upper complex half-plane, has
slope equal to 1 (resp. to -1), and cuts the real axis between 0 and 1.
Setting in (8.9) $s = \hf + it$ and using the property (8.4) it
follows that
$$
Z(t) \> = \> 2\,\Re{\rm e}\Bigl({\rm e}^{-i\theta(t)}
\int\displaylimits_{0\searrow 1}{{\rm e}^{-i\pi z^2}z^{-1/2+it}\over
{\rm e}^{i\pi z} - {\rm e}^{-i\pi z}}\d z\Bigr).
$$\noindent
As $\Re{\rm e}\,(iw) = -\Im{\rm m}\,w$, this can be conveniently written as
$$
Z(t) \> = \> \Im{\rm m}\Bigl({\rm e}^{-i\theta(t)}
\int\displaylimits_{0\searrow 1}{\rm e}^{-i\pi z^2}z^{-1/2+it}\,
{\d z\over\sin(\pi z)}\Bigr).
\leqno{(8.10)}
$$\noindent
Since
$$
\vert\sin z\vert^2 = \sin^2(\Re{\rm e}\,z) + \sinh^2(\Im{\rm m}\,z),\>
z^{-1/2+it} = \vert z\vert^{-1/2+it}{\rm e}^{-{i\over 2}\arg z - t\arg z},
$$\noindent
and for $z = \eta + u{\rm e}^{3\pi i/4}$, $u$ real, $0 < \eta < 1$ we have
$$
\vert {\rm e}^{-i\pi z^2}\vert \> = \> {\rm e}^{-\pi u^2 + \eta\sqrt{2}\pi u},
$$\noindent
it follows that the contribution of the portion of the integral in (8.10)
for which $\vert z\vert \ge \log t$ is $\ll {\rm e}^{-\log^2t}$. Hence
$$
Z(t) \> = \> \Im{\rm m}\Bigl({\rm e}^{-i\theta(t)}
\int\displaylimits_{0\searrow 1,\vert z\vert < \log t}
{\rm e}^{-i\pi z^2}z^{-1/2+it}\,
{\d z\over\sin(\pi z)}\Bigr) + O({\rm e}^{-\log^2t}). \leqno{(8.11)}
$$ \noindent
From the decay property (6.4) it follows that
$$
M_{Z,f}(t) = \int_{-\infty}^\infty Z(t + x)f\bigl({x\over G}\bigr)\d x
= \int\displaylimits_{-\log^{2\alpha -1}t}^{\log^{2\alpha -1}t} 
Z(t + x)f\bigl({x\over G}\bigr)\d x
+   O({\rm e}^{-c\log^2t}),\leqno{(8.12)}
$$ \noindent
where $c$ denotes positive, absolute constants which may not be the same
ones at each occurrence. Thus from (8.10) and (8.12) we obtain that
$$
M_{Z,f}(t) =             \Im{\rm m}\Bigl(
\int\displaylimits_{0\searrow 1,\vert z\vert < \log t}
{{\rm e}^{-i\pi z^2}z^{-1/2+it}
\over\sin(\pi z)}\int\displaylimits_{-\log^{2\alpha -1}t}
^{\log^{2\alpha -1}t}
{\rm e}^{-i\theta(t+x)}z^{ix}f\bigl({x\over G}\bigr)\d x\d z\Bigr)
+ O({\rm e}^{-c\log^2t}).\leqno{(8.13)}
$$
\noindent
By using Taylor's formula we have
$$
\theta(t + x)\>=\>\theta(t) + {x\over 2}\log{t\over 2\pi} + x\Delta'(t) 
+ R(t,x) \leqno{(8.14)}
$$ \noindent
with $\Delta'(t) \ll t^{-2}$ and
$$
R(t,x) = \sum_{n=2}^\infty \biggl({(-1)^n\over 2n(n - 1)t^{n-1}}
+ {\Delta^{(n)}(t)\over n!}\biggr)x^n.
$$
\noindent
Now we put
$$
{\rm e}^{-iR(t,x)} = 1 + S(t,x), \leqno{(8.15)}
$$\noindent
say, and use (8.5), (8.6), (8.8) and (8.14). We obtain
$$
S(t,x) = \sum_{k=1}^\infty{(-i)^kR^k(t,x)\over k!}
= \sum_{n=2}^\infty g_n(t)x^n, \leqno{(8.16)}
$$ \noindent
where each $g_n(t) \;(\in C^\infty (0,\infty))$ has an asymptotic
expansion of the form
$$
g_n(t) \>\sim\> \sum_{k=0}^\infty d_{n,k}t^{-k-[(n+1)/2]} \quad (t \to\infty)
\leqno{(8.17)}
$$\noindent
with suitable constants $d_{n,k}$. From (8.13)-(8.15) we have
$$
M_{Z,f}(t) \> =\> \Im{\rm m}\,\Bigl(I_1 + I_2\Bigr) 
+ O({\rm e}^{-c\log^2t}), \leqno{(8.18)}
$$\noindent
where
$$
I_1 := 
\int\displaylimits_{0\searrow 1,\vert z\vert < \log t}
{{\rm e}^{-i\pi z^2}z^{-1/2+it}
\over\sin(\pi z)}\int\displaylimits_{-\log^{2\alpha -1}t}
^{\log^{2\alpha -1}t}
{\rm e}^{-i\theta(t)-{ix\over 2}\log{t\over 2\pi}-ix\Delta'(t)}
z^{ix}f\bigl({x\over G}\bigr)\d x\d z, \leqno{(8.19)}
$$
$$
I_2 := 
\int\displaylimits_{0\searrow 1,\vert z\vert < \log t}
{{\rm e}^{-i\pi z^2}z^{-1/2+it}
\over\sin(\pi z)}\int\displaylimits_{-\log^{2\alpha -1}t}
^{\log^{2\alpha -1}t}
{\rm e}^{-i\theta(t)-{ix\over 2}\log{t\over 2\pi}-ix\Delta'(t)}
z^{ix}f\bigl({x\over G}\bigr)S(t,x)\d x\d z. \leqno{(8.20)}
$$            \noindent
In $I_1$ we write $z = {\rm e}^{i\arg z}\vert z\vert$, which gives
$$
I_1 := 
\int\displaylimits_{0\searrow 1,\vert z\vert < \log t}
{\rm e}^{-\theta(t)}{{\rm e}^{-i\pi z^2}z^{-1/2+it}
\over\sin(\pi z)}h(z)\d z + O({\rm e}^{-c\log^2t}),
$$
\noindent
where
$$ \eqalign{
h(z) &:= \int_{-\infty}^\infty {\rm e}^{-x\arg z}f\bigl({x\over G}\bigr)
\,\exp\biggl(ix\log{\vert z\vert {\rm e}^{-\Delta'(t)}
\over\sqrt{t/2\pi}}\biggr)\d x           \cr&
=G \int_{-\infty}^\infty {\rm e}^{-Gy\arg z}f(y)
\,\exp\biggl(iyG\log{\vert z\vert {\rm e}^{-\Delta'(t)}
\over\sqrt{t/2\pi}}\biggr)\d y\cr&
= G\sum_{n=0}^\infty {(-G\arg z)^n\over n!}
 \int_{-\infty}^\infty y^nf(y)
\,\exp\biggl(2i\pi y\cdot {G\over 2\pi}\log{\vert z\vert
{\rm e}^{-\Delta'(t)}
\over\sqrt{t/2\pi}}\biggr)\d y.\cr}
$$  \noindent
Change of summation and integration is justified by
absolute convergence, since $\vert z\vert < \log t,
G = \delta/\log(T/2\pi)$, $\Delta'(t) \ll t^{-2}, t \sim T$,
and $f(x)$ satisfies (6.4). But
$$
 \int_{-\infty}^\infty f(y)
\,\exp\biggl(2i\pi y\cdot {G\over 2\pi}\log{\vert z\vert {\rm e}^{-\Delta'(t)}
\over\sqrt{t/2\pi}}\biggr)\d y
= \widehat f
\biggl({G\over 2\pi}\log{\vert z\vert {\rm e}^{-\Delta'(t)}
\over\sqrt{t/2\pi}}\biggr) = 1
$$ \noindent
for $\delta < 2\pi (b - a)$, since
$$\Big\arrowvert {G\over 2\pi}\log{\vert z\vert {\rm e}^{-\Delta'(t)}
\over\sqrt{t/2\pi}} \Big\arrowvert =
{\delta\over 2\pi\log\bigl({T\over 2 \pi}\bigr)}
\Bigl(\log\sqrt{t\over 2\pi} - \log \vert z\vert + \Delta'(t)\Bigr)
= ({\delta\over 4\pi} + o(1)) < {\delta\over 2\pi} < b - a,
$$ \noindent
and $\widehat f(x) = 1$ for $\vert x \vert < b - a$. Moreover
for $n \ge 1$ and $\vert x \vert < b - a$ we have
$$
\widehat f^{(n)}(x) = (2\pi i)^n\int_{-\infty}^\infty
y^ne^{2\pi ixy}f(y)\d y = 0,
$$\noindent
hence it follows that
$$
\int_{-\infty}^\infty y^nf(y)
\,\exp\biggl(2i\pi y\cdot {G\over 2\pi}\log{\vert z\vert
{\rm e}^{-\Delta'(t)}
\over\sqrt{t/2\pi}}\biggr)\d y = 0 \quad (n \ge 1, 1 < 2\pi (b - a)).
$$\noindent
Thus we obtain
$$
I_1 = G
\int\displaylimits_{0\searrow 1}
{{\rm e}^{-\theta(t)}{\rm e}^{-i\pi z^2}z^{-1/2+it}
\over\sin(\pi z)}\d z + O({\rm e}^{-c\log^2t}). \leqno{(8.21)}
$$
\noindent
Similarly from (8.16) and (8.20) we have
$$
I_2 = 
\int\displaylimits_{0\searrow 1,\vert z\vert < \log t}
{{\rm e}^{-i\pi z^2}z^{-1/2+it}
\over\sin(\pi z)}\int\displaylimits_{-\log^{2\alpha -1}t}
^{\log^{2\alpha -1}t}
{\rm e}^{-i\theta(t)-{ix\over 2}\log{t\over 2\pi}-ix\Delta'(t)}
z^{ix}f\bigl({x\over G}\bigr)\sum_{n=2}^\infty x^ng_n(t)\d x\d z
$$
$$
= \int\displaylimits_{0\searrow 1,\vert z\vert < \log t}
{{\rm e}^{-i\pi z^2}z^{-1/2+it}
\over\sin(\pi z)}\int\displaylimits_{-\log^{2\alpha -1}t}
^{\log^{2\alpha -1}t}
\biggl(\sum_{n=1}^N{P_n(x)\over t^n} + 
 O\bigl({1 + x^{N+2}\over t^{N+1}}\bigr)\biggr)
{\rm e}^{-i\theta(t)-{ix\over 2}\log{t\over 2\pi}-ix\Delta'(t)}
z^{ix}f\bigl({x\over G}\bigr)\d x\d z
$$
$$
= \sum_{n=1}^N e^{-i\theta(t)}t^{-n}
\int\displaylimits_{0\searrow 1,\vert z\vert < \log t}
{{\rm e}^{-i\pi z^2}z^{-1/2+it}
\over\sin(\pi z)}\int\displaylimits_{-\log^{2\alpha -1}t}
^{\log^{2\alpha -1}t}P_n(x)
{\rm e}^{-{ix\over 2}\log{t\over 2\pi}-ix\Delta'(t)}
z^{ix}f\bigl({x\over G}\bigr)\d x\d z
$$
$$
+\; O\Biggl({1\over t^{N+1}}
\int\displaylimits_{-\log^{2\alpha -1}t}^{\log^{2\alpha -1}t}
(1 + x^{N+2})\vert f\bigl({x\over G}\bigr)\vert\,
\Big\arrowvert 
\int\displaylimits_{0\searrow 1}{{\rm e}^{-i\pi z^2}z^{-1/2+it+ix}
\over\sin(\pi z)}\,dz\Big\arrowvert \d x\,\Biggr) + O({\rm e}^{-c\log^2t}),
$$     \noindent
where each $P_n(x)$ is a polynomial in $x$ of degree $n (\ge 2)$.
The integral over $z$ in the error term is similar to the one in (8.10).
Hence by the residue theorem we have ($Q = [\sqrt{t/2\pi}\,]$)
$$
\int\displaylimits_{0\searrow 1}{{\rm e}^{-i\pi z^2}z^{-1/2+it+ix}
\over\sin(\pi z)}\d z =
2\pi i\sum_{n=1}^Q \mathop{\rm Res}\limits_{z=n} 
{{\rm e}^{-i\pi z^2}z^{-1/2+it+ix} \over\sin(\pi z)}  + 
\int\displaylimits_{Q\searrow Q+1}{{\rm e}^{-i\pi z^2}z^{-1/2+it+ix}
\over\sin(\pi z)}\d z,
$$ \noindent
similarly as in the derivation of the Riemann-Siegel formula. It follows that
$$
\Big\arrowvert 
\int\displaylimits_{0\searrow 1}{{\rm e}^{-i\pi z^2}z^{-1/2+it+ix}
\over\sin(\pi z)}\,dz\Big\arrowvert
\;\ll\;  t^{1/4}.
$$ \noindent
Thus analogously as in the case of $I_1$ we find that, for $n \ge 1$,
$$
\int\displaylimits_{-\log^{2\alpha -1}t}^{\log^{2\alpha -1}t}
P_n(x){\rm e}^{-i\theta(t)-{ix\over 2}\log{t\over 2\pi}-ix\Delta'(t)}
z^{ix}f\bigl({x\over G}\bigr)\d x
= \int_{-\infty}^\infty P_n(x)\cdots\d x 
+  O({\rm e}^{-c\log^2t}) = O({\rm e}^{-c\log^2t}).
$$\noindent
Hence it follows that, for any fixed integer $N \ge 1$,
$$
I_2 \>\ll_N \>T^{-N}. \leqno{(8.22)}
$$\noindent
Theorem 1 now follows from (8.10) and (8.18)-(8.22), since clearly
it suffices to assume that $N$ is an integer. One can generalize
Theorem 1 to derivatives of $M_{Z,f}(t)$.
\medskip

Theorem 1 shows that $Z(t)$ and $M_{Z,f}(t)/G$ differ only by $O(T^{-N})$,
for any fixed $N \ge 0$, which is a very small quantity. 
This certainly supports the belief that, for this particular subclass of
functions $f(x)$, (7.8)--(7.11) will be true, but {\it proving}
it may be very hard. On the other hand, nothing precludes the 
possibility that the error
term in Theorem 1, although it is quite small, represents a function 
possessing many small ``spikes" (like $t^{-N}\sin(t^{N+2})$, say).
These spikes could introduce many new zeros, thus violating (7.8)--(7.11). 
Therefore it remains an open question to investigate the distribution of
zeros of $M_{Z,f}(t)$ of Theorem 1, and to see to whether there is a 
possibility  that Theorem 1 can be used in settling the truth of the RH.

\bigskip

\centerline{\dd 9. Convolution functions and the RH}

\bigskip

In this section we shall discuss the possibility to use 
convolution functions to  disprove the RH, of course assuming that
it is false. Let us denote by $T_\alpha^\beta$ the subclass of
$S_\alpha^\beta$  with $\alpha > 1$ consisting of 
functions $f(x)$, which are not identically equal to zero, and which satisfy 
$\int_{-\infty}^\infty f(x)\d x \; >\;0.$
It is clear that $T_\alpha^\beta$ is non-empty. Our choice for $G$
will be the same one as in (8.1),
so that for suitable $\delta$ we shall have
$$
\widehat{f}\Bigl({G\over 4\pi}\log\bigl({T\over 2\pi}\bigr)\Bigr) =
\widehat{f}\bigl({\delta\over 4\pi}\bigr) \gg 1. \leqno{(9.1)}
$$\noindent
In fact by continuity (9.1) will hold for $\vert\delta\vert \le C_1$,
where $C_1 > 0$ is a suitable constant depending only on $f$, since
if  $f(x) \in T_\alpha^\beta$, then we have
$$
\widehat{f}(0)\;=\;\int_{-\infty}^\infty f(x)\d x \; >\;0.\leqno{(9.2)}
$$\noindent
 Moreover if $f(x) \in S_\alpha^0$, then  $\widehat{f}(x)
\in S_0^\alpha$ and thus it is of bounded support, and consequently
$G \ll 1/\log T$ must hold if the bound in (9.1) is to be satisfied.
This choice of $f(x)$  turns out to be better suited for our purposes
than the choice made in Section 8, which perhaps would seem more natural in
view of Theorem 1. The reason for this is that, if $f(x) \in S_\alpha^\beta$
with parameters $A$ and $B$, then $\widehat{f}(x) \in S_\beta^\alpha$ with
parameters $B + \varepsilon$ and $A + \varepsilon$, respectively (see [9]).
But for $f(x)$ as in Section 8 we have $A = a + b$, thus for $\widehat{f}(x)$
we would have (in [9] $\widehat{f}$ is defined without the factor $2\pi$, 
which would only change the scaling factors) $B = a + b + \varepsilon$, and
this value of $B$ would eventually turn out to be too large for our 
applications. In the present approach we have more flexibility, since only
(9.2) is needed. Note that $f(x)$ is not necessarily nonnegative.
\smallskip
Now observe that if we replace $f(x)$ by
$f_1(x) := f(Dx)$ for a given $D > 0$, then obviously $f_1(x) \in
S_\alpha^0$, and moreover uniformly for $q \ge 0$ we have 
$$
f_1^{(q)}(x) = D^qf^{(q)}(Dx) \ll (BD)^q\exp(-aD^{1/\alpha}
\vert x\vert^{1/\alpha}).
$$
\noindent
In other words the constant $B$ in (6.3)
or (6.4) is replaced by $BD$. 
Take now $D = \eta /B$, where $\eta > 0$ is an arbitrary, but fixed
number, and write $f$ for $Df_1$. If the RH holds, 
then from (4.5), (6.4) and (6.5) and we have, for $k$ given by (7.12),
$$
M_{Z,f}^{(k)}(t) \ll {\Bigl( {\eta\over G}\Bigr)}^k\exp\Bigl({B_1
\log t\over \log\log t}\Bigr)\leqno{(9.3)}
$$ \noindent
with a suitable constant $B_1 > 0$.

We shall assume now that the RH holds and that (7.8), (7.11) 
hold for some  $f(x) \in T_\alpha^0$
(for which (9.3) holds, which is implied by the RH),  
and we shall obtain a contradiction. To this end
let  $\,U := T^{1/2+\varepsilon}\,$, so that we may apply (7.10) or
(7.11), $V = T^{1/4}/\log T, L = \log^{1/2+\varepsilon}T$. We shall
consider the mean value of $\vert M_{Z,f}(t)\vert$ over $\,[T - U, T +
U]\,$ in order to show that, on the average,  $\vert M_{Z,f}(t)\vert$
is not too small. We have
$$\eqalign{
I &: = \int\displaylimits_{T-U}^{T+U}\vert M_{Z,f}(t)\vert\d t \ge
\sum_{n=1}^N \>\int\displaylimits_{T_n-VL}^{T_n+VL}\vert M_{Z,f}(t)\vert\d t
\cr&
\ge \; \sum_{n=1}^N\> \int\displaylimits_{T_n-VL}^{T_n+VL}
\vert M_{Z,f}(t)\vert
\exp\Bigl(-(T_n - t)^2V^{-2}\Bigr)\d t,\cr}\leqno{(9.4)}
$$ \noindent
where $T_n = T - U + (2n - 1)VL$, and $N$ is the largest integer
for which $T_N + VL \le T + U$, hence $N \gg UV^{-1}L^{-1}$. We use
Lemma 1 to bound from below each integral over $\,[T_n - VL, T_n + VL]\,$.
It follows that
$$
I \gg GV\sum_{n=1}^N\Bigl(\Big\arrowvert\widehat{f}\bigl({G\over
4\pi}\log\bigl({T_n\over 2\pi}\bigr)\bigr)\Big\arrowvert + O(T^{-1/4})\Bigr)
\gg  GV(N + O(NT^{-1/4})) \gg GUL^{-1} \leqno{(9.5)}
$$ \noindent
for sufficiently small $\delta$, since  (9.1) holds and
$$
\widehat{f}\Bigl({G\over 4\pi}\log\bigl({T_n\over 2\pi}\bigr)\Bigr) =
\widehat{f}\Bigl({\delta\over4\pi}\cdot{\log\bigl({T_n\over
2\pi}\bigr)\over\log\bigl({T\over 2\pi}\bigr)}\Bigr) =
\widehat{f}\Bigl({\delta\over 4\pi} + O\bigl({U\over
T}\bigr)\Bigr).
$$

\medskip
We have assumed that (7.11) holds, but this implies that (7.10) holds
also. Namely it holds unconditionally with $S(t)$ in place of
$S_M(t)$.  Thus for any fixed integer $m \ge 1$ we have
$$
\int_T^{T+U}\bigl(S_M(t + H) - S_M(t - H)\bigr)^{2m}\d t
$$
$$
\ll \int_{T+H}^{T+H+U}\bigl(S_M(t) - S(t)\bigr)^{2m}\d t
+ \int_T^{T+U}\bigl(S(t + H) - S(t - H)\bigr)^{2m}\d t
+ \int_{T-H}^{T-H+U}\bigl(S(t) - S_M(t)\bigr)^{2m}\d t
$$
$$
\ll U(\log\log T)^m \qquad (T^a < U \le T,  {1\over 2} < a \le 1).
$$

Let $\cal D$ be the subset of $\,[T - U, T + U]\,$ where 
$$
\vert S_M(t + H) - S_M(t - H)\vert  \;\le\; \log^{1/2}T\leqno{(9.6)}
$$ \noindent
fails. The bound (7.10) implies that
$$
m({\cal D}) \ll U\log^{-C}T \leqno{(9.7)}
$$ \noindent
for any fixed $C > 0$. If we take $C = 10$ in (9.7) and use the
Cauchy-Schwarz inequality for integrals we shall have
$$
\int_{\cal D}\vert M_{Z,f}(t)\vert\d t  \le (m({\cal
D}))^{1/2}\Bigl(\int_{T-U}^{T+U}M_{Z,f}^2(t)\d t\Bigr)^{1/2}
\ll GU\log^{-4}T, \leqno{(9.8)}
$$ \noindent
since
$$ \eqalign{&
\int_{T-U}^{T+U}M_{Z,f}^2(t)\d t
\le \int_{T-U}^{T+U}\int_{-\infty}^\infty Z^2(t + x)|f({x\over G})|\d x 
\int_{-\infty}^\infty |f({y\over G})|\d y\,\d t    \cr&
\ll G\int_{T-U}^{T+U}\int\displaylimits_{-G\log^{2\alpha -1}T}
^{G\log^{2\alpha -1}T}
\vert\zeta(\hf + it + ix)\vert^2|f({x\over G})|\d x\,\d t + G\cr&
\ll G\int\displaylimits_{-G\log^{2\alpha -1}T}^{G\log^{2\alpha -1}T}
\Bigl(\int\displaylimits_{T-2U}^{T+2U}
\vert\zeta(\hf + iu)\vert^2\d u\Bigr)|f({x\over G})|\d x + G
\ll G^2U\log T.                                        \cr}
$$ \medskip \noindent
The last bound easily follows from mean square results on
$\vert\zeta({1\over 2} + it)\vert$ (see [16]) with the
choice $U = T^{{1\over 2}+\varepsilon}$. Therefore (9.4) and
(9.8) yield
$$
GUL^{-1} \ll \int\displaylimits_{{\cal D}'}\vert M_{Z,f}(t)\vert\,dt, 
\leqno{(9.9)}
$$ \medskip\noindent
where ${\cal D}' = [T - U, T + U] \; \backslash \;{\cal D}$, hence
in (9.9) integration is over $t$ for which (9.6) holds. If $t \in {\cal D'}$,
$\gamma$ denotes the zeros of $M_{Z,f}(t)$, then from (7.7) and (7.12)
we obtain (recall that $\log_rt = \log(\log_{r-1}t)$)
$$
k = k(t,T) = {H\over \pi}\log({T\over 2\pi})\cdot\{1 + O((\log
T)^{\varepsilon-{1\over 2}})\}, \; \log k = \log H - \log\pi
+ \log_2({T\over2\pi}) + O((\log T)^{\varepsilon-1/2})
\leqno{(9.10)}
$$ \noindent
for any given $\varepsilon > 0$. To bound $M_{Z,f}(t)$ we use (7.6),
with $k$ given by (9.10), $\tau = \tau(t,k)$, (9.3) and
$$
k! = \exp(k\log k - k + O(\log k)).
$$ \noindent
We obtain, denoting by $B_j$ positive absolute constants,
$$
GUL^{-1} \ll \int\displaylimits_{{\cal D}'}\prod_{\vert\gamma - t\vert\le H}
\vert\gamma - t\vert{\vert M_{Z,f}^{(k)}(\tau )\vert\over k!}\d t
\leqno{(9.11)}
$$
$$
\ll \exp\Bigl({B_2\log T\over \log_2T}\Bigr)\int\displaylimits_{{\cal D}'}
\exp\{k(\log{\eta \over G} - \log k + 1)\}\prod_{\vert\gamma - t\vert\le H}
\vert\gamma - t\vert\d t
$$
$$
\ll \exp\Bigl({B_3\log T\over \log_2T} + {H\over \pi}\log({T\over
2\pi})\cdot\bigl(\log{\eta \over \delta} +
\log_2({T\over 2\pi}) - \log H + \log\pi
- \log_2({T\over 2\pi}) + 1\bigr)\Bigr)\,
\int\displaylimits_{{\cal D}'}\prod_{\vert\gamma - t\vert\le H}
\vert\gamma - t\vert\d t.
$$ \medskip\noindent
It was in evaluating $k\log k$ that we needed (9.10), since only the
bound (7.9) would not suffice (one would actually need the bound
$S_M(T) \ll \log T/(\log _2T)^2$). If the product under the last
integral is bounded by (7.13), we obtain

\medskip
$$
GUL^{-1} \ll U\exp\Bigl({H\over \pi}\log({T\over 2\pi})\cdot (\log{\eta \over \delta}
+ B_4)\Bigr),
$$ \noindent
and thus for $T \ge T_0$
$$
1 \le \exp\Bigl({H\over \pi}\log({T\over 2\pi})\cdot (\log{\eta \over \delta}
+ B_5)\Bigr). \leqno{(9.12)}
$$\noindent
Now we choose e.g.,
$$
\eta = \delta^2,\;\delta \;=\;\min\,(C_1,\,{\rm e}^{-2B_5}),
$$\noindent
where $C_1$ is the constant for which (9.1) holds if 
$\vert\delta\vert \le C_1$, so that (9.12) gives
$$
1 \le \exp\Bigl({-B_5H\over \pi}\log({T\over 2\pi})\Bigr),
$$\noindent
which is a contradiction for $T \ge T_1$.  Thus we have proved the following

\bigskip \medskip
{\bf THEOREM 2}. {\it If} (7.8) {\it and} (7.11) {\it hold for suitable
$f(x) \in T_\alpha^0$ with $G$ given by} (8.1), {\it then the Riemann
hypothesis is false}.

\bigskip
Theorem 2 is  similar to the result proved also in [23]. 
Actually the method of  proof of Theorem 2 gives more than the
assertion of the theorem. Namely it shows that, under the above
hypotheses, (4.5) cannot hold for any fixed $A > 0$.
Perhaps it should be
mentioned that (7.11) is not the only condition which would lead to the
disproof of the RH. It would be enough to assume, under the RH, that one
had (7.8)--(7.10) for a suitable $f(x)$, or
$$
N_M(t) \;=\; N(t) + O\Bigl({\log T\over (\log\log T)^2}\Bigr)
\leqno{(9.18)}
$$ \noindent
for $t \in [T - U, T + U]$ with a suitable $U (= T^{1/2+\varepsilon}$,
but smaller values are possible), to derive a contradiction. The main
drawback of this approach is the necessity to impose conditions like
(7.8)--(7.10) which can be, for all we know, equally difficult to settle 
as the assertions which we originaly set out to prove (or disprove). 
For this reason our results can only be conditional. Even if the RH
is false it appears plausible that, as $T \to \infty$, $N_0(T) = (1 +
o(1))N(T)$. In other words, regardless of the truth of the RH, almost
all complex zeros of $\zeta(s)$ should lie on the critical line. 
This is the conjecture that the author certainly believes in. No plausible
conjectures seem to exist (if the RH is false) regarding the order of
$N(T) - N_0(T)$.

\bigskip        \bigskip
\vfill
\eject

\topglue2cm
\centerline{\dd REFERENCES}

\bigskip\bigskip

\item{ [1]} R.J. Anderson and H.M. Stark, {\it Oscillation theorems},
in LNM's {\bf 899}, Springer-Verlag, Berlin-Heidelberg-New York,
1981, 79-106.

\item{ [2]} F.V. Atkinson, {\it The mean value of the Riemann
zeta-function}, Acta Math. {\bf 81}(1949), 353-376.

\item{ [3]} R. Balasubramanian, {\it On the frequency of Titchmarsh's
phenomenon for $\zeta(s)$ IV}, Hardy-Ramanujan J.
{\bf 9}(1986), 1-10.

\item{ [4]} R. Balasubramanian and K. Ramachandra, {\it On the frequency
of Titchmarsh's phenomenon for $\zeta(s)$ III},
 Proc. Indian Acad. Sci. Section A {\bf 86}(1977), 341-351.

\item{ [5]} E. Bombieri and D. Hejhal, {\it Sur les z\'eros des
fonctions zeta d'Epstein}, Comptes Rendus Acad. Sci. Paris
{\bf 304}(1987), 213-217.

\item{ [6]} H. Davenport and H. Heilbronn, {\it On the zeros of certain
Dirichlet series I,II}, J. London Math. Soc.
 {\bf 11}(1936), 181-185 and ibid. 307-312.

\item{ [7]} H.M. Edwards, {\it Riemann's zeta-function}, Academic Press,
New York-London, 1974.

\item{ [8]} A. Fujii, {\it On the distribution of the zeros of the
Riemann zeta-function in short intervals}, Bull. Amer. Math. Soc.
{\bf 8}1(1975), 139-142.

\item{ [9]} I.M. Gel'fand and G.E. Shilov, {\it Generalized functions
(vol. 2)}, Academic Press, New York-London, 1968.

\item{ [10]} A.O. Gel'fond, {\it The calculus of finite differences
(Russian)}, Nauka, Moscow, 1967.

\item{ [11]} J.L. Hafner and A. Ivi\'c, {\it On the mean square of the
Riemann zeta-function on the critical line}, J. Number Theory {\bf 32}(1989),
151-191.

\item{ [12]  }   J.L. Hafner and A. Ivi\'c, {\it On some mean value results
for the Riemann zeta-function}, in Proc. International  Number Theory Conf.
Qu\'ebec 1987, Walter de Gruyter and Co., Berlin-New York, 348-358.

\item{ [13]} M.N. Huxley, {\it Exponential sums and the Riemann zeta-function
IV}, Proceedings London Math. Soc. (3){\bf 66}(1993), 1-40.

\item{ [14]} M.N. Huxley, {\it A note on exponential sums with a
difference}, Bulletin London Math. Soc. {\bf 29}(1994), 325-327.

\item{ [15]} A. Ivi\'c, {\it Large values of the error term in the
divisor problem}, Invent. Math. {\bf 71}(1983), 513-520.

\item{ [16]} A. Ivi\'c, {\it The Riemann zeta-function}, John Wiley
and Sons, New York, 1985 (2nd ed. Dover, 2003).

\item{ [17]} A. Ivi\'c, {\it On consecutive zeros of the Riemann
zeta-function on the critical line}, S\'eminaire de Th\'eorie des
Nombres, Universit\'e de Bordeaux 1986/87, Expos\'e no. {\bf 29}, 14 pp.

\item{ [18]} A. Ivi\'c, {\it On a problem connected with zeros of
$\zeta(s)$ on the critical line}, Monatshefte Math. {\bf 104}(1987), 17-27.

\item{ [19]} A. Ivi\'c, {\it Large values of certain number-theoretic
error terms}, Acta Arithmetica {\bf 56}(1990), 135-159.

\item{ [20]} A. Ivi\'c, {\it Mean values of the Riemann zeta-function},
LN's {\bf 82}, Tata Institute of Fundamental Research, Bombay,
1991 (distr. by Springer Verlag, Berlin etc.).

\item{ [21]} A. Ivi\'c, {\it On a class of convolution functions connected
with $\zeta(s)$}, Bulletin CIX Acad. Serbe des Sciences
 et des Arts, Classe des Sciences math\'ematiques et naturelles,
Math. No {\bf 20}(1995), 29-50.

\item{ [22]} A. Ivi\'c, {\it On the fourth moment of the Riemann
zeta-function}, Publications Inst. Math.
(Belgrade) {\bf 57(71)}(1995), 101-110.

\item{ [23]} A. Ivi\'c, {\it On the distribution of zeros of a class
of convolution functions}, 
        Bulletin CXI Acad. Serbe des Sciences et des
 Arts, Classe des Sciences math\'ematiques
et naturelles,  Sciences math\'ematiques No. {\bf21}(1996), 61-71.

\item{ [24]} A. Ivi\'c, {\it On the ternary additive divisor problem and
 the sixth moment of the
zeta-function}, in ``Sieve Methods, Exponential Sums, and
their Applications in
Number Theory" (eds. G.R.H. Greaves, G. Harman, M.N. Huxley), Cambridge
University Press (Cambridge, UK), 1996, 205-243.

\item{ [25]} A. Ivi\'c and M. Jutila, {\it Gaps between consecutive
zeros of the Riemann zeta-function}, Monatshefte Math.
 {\bf 105}(1988), 59-73.

\item{ [26]} A. Ivi\'c and Y. Motohashi, {\it A note on the mean value of
the zeta and L-functions VII}, Proc. Japan Acad.
 Ser. A {\bf 66}(1990), 150-152.

\item{ [27]} A. Ivi\'c and Y. Motohashi, {\it The mean square of the
error term for the fourth moment of the zeta-function},
 Proc. London Math. Soc. (3){\bf 66}(1994), 309-329.

\item{ [28]} A. Ivi\'c and Y. Motohashi, {\it The fourth moment of the
Riemann zeta-function}, Journal Number Theory {\bf 51}(1995), 16-45.

\item{ [29]} A. Ivi\'c and H.J.J. te Riele, {\it On the zeros of the error
term for the mean square of $\;\vert\zeta({1\over 2}  + it)\vert\;$}, Math.
Comp. {\bf 56} No {\bf 193}(1991), 303-328.

\item{ [30]} M. Jutila, {\it On the value distribution of the zeta-function
on the critical line}, Bull. London Math. Soc.
{\bf 15}(1983), 513-518.

\item{ [31]} A.A. Karatsuba, {\it On the zeros of the Davenport-Heilbronn
function lying on the critical line (Russian)}, Izv. Akad. Nauk SSSR
ser. mat. {\bf 54} no. 2 (1990), 303-315.

\item{ [32]} A.A. Karatsuba and S.M. Voronin, {\it The Riemann
zeta-function}, Walter de Gruyter, Berlin--New York, 1992.

\item{ [33]} N.V. Kuznetsov, {\it Sums of Kloosterman sums and the
eighth moment of the Riemann zeta-function}, Papers presented at the
Ramanujan Colloquium, Bombay 1989, publ. for Tata Institute (Bombay)
by Oxford University Press, Oxford, 1989, pp. 57-117.

\item{ [34]} A.A. Lavrik, {\it Uniform approximations and zeros of
derivatives of Hardy's Z-function in short intervals (Russian)},
Analysis Mathem. {\bf 17}(1991), 257-259.

\item{ [35]} D.H. Lehmer, {\it On the roots of the Riemann zeta
function}, Acta Math. {\bf 95}(1956), 291-298.

\item{ [36]} D.H. Lehmer, {\it Extended computation of the Riemann
zeta-function}, Mathematika {\bf 3}(1956), 102-108.

\item{ [37]} J.E. Littlewood, {\it Sur la distribution des nombres
premiers}, Comptes rendus Acad\'emie Sci. (Paris) {\bf 158}(1914),
 1869-1872.

\item{ [38]} J. van de Lune, H.J.J. te Riele and D.T. Winter,
{\it On the zeros of the Riemann zeta-function in the critical strip IV},
Math. Comp. {\bf 46}(1987), 273-308.

\item{ [39]}  H.L. Montgomery, {\it Extreme values of the Riemann
zeta-function}, Comment. Math. Helv. {\bf 52}(1977), 511-518.

\item{ [40]} Y. Motohashi, {\it The fourth power mean of the Riemann
zeta-function}, in "Proceedings of the Amalfi
 Conference on Analytic
Number Theory 1989", eds. E. Bombieri et al., Universit\`a di Salerno,
Salerno, 1992, 325-344.

\item{ [41]} Y. Motohashi,  {\it  An explicit formula for the fourth
power mean of the Riemann zeta-function}, Acta Math.
{\bf 170}(1993), 181-220.

\item{ [42]} Y. Motohashi, {\it A relation between the Riemann zeta-function
and the hyperbolic Laplacian}, Ann. Sc.
 Norm. Sup. Pisa, Cl. Sci. IV
ser. {\bf 22}(1995), 299-313.

\item{ [43]} Y. Motohashi, {\it The Riemann zeta-function and the
non-Euclidean Laplacian}, Sugaku Expositions, AMS
{\bf 8}(1995), 59-87.

\item{ [44]} A.M. Odlyzko, {\it On the distribution of spacings between
the zeros of the zeta-function}, Math. Comp.
  {\bf 48}(1987), 273-308.

\item{ [45]} A.M. Odlyzko, {\it Analytic computations in number theory},
Proc. Symposia in Applied Math. {\bf 48}(1994), 451-463.

\item{ [46]} A.M. Odlyzko, {\it The} $10^{20}$-{\it th zero of the Riemann
zeta-function and 175 million of its neighbors}, to appear.

\item{ [47]} A.M. Odlyzko and H.J.J. te Riele, {\it Disproof of the
Mertens conjecture}, J. reine angew. Math. {\bf 357}(1985), 138-160.

\item{ [48]} K. Ramachandra, {\it Progress towards a conjecture on the
mean value of Titchmarsh series}, in "Recent
 Progress in Analytic Number
Theory", symposium Durham 1979 (Vol. 1), Academic Press, London, 1981,
303-318.

\item{ [49]} K. Ramachandra, {\it On the mean-value and omega-theorems
for the Riemann zeta-function}, LNs 85, Tata
Institute of Fundamental
Research, Bombay, 1995 (distr. by Springer Verlag, Berlin etc.).

\item{ [50]} H.J.J. te Riele, {\it On the sign of the difference} $\pi(x) -
{\rm li}\,x$, Math. Comp. {\bf 48}(1987), 323-328.

\item{ [51]} H.J.J. te Riele and J. van de Lune, {\it Computational  number
theory  at CWI in 1970-1994}, CWI Quarterly {\bf 7(4)} (1994), 285-335.

\item{ [52]} B. Riemann, {\it \"Uber die Anzahl der Primzahlen unter
einer gegebener Gr\"osse}, Monats. Preuss. Akad. Wiss. (1859-1860),
671-680.

\item{ [53]} A. Selberg, {\it On the zeros of  Riemann's zeta-function},
Skr. Norske Vid. Akad. Oslo {\bf 10}(1942), 1-59.

\item{ [54]} A. Selberg, {\it Contributions to the theory of the Riemann
zeta-function}, Arch. Math. Naturvid. {\bf 48}(1946) No. 5, 89-155.

\item{ [55]} A. Selberg, {\it Harmonic analysis and discontinuous
groups in weakly symmetric spaces with applications
to Dirichlet series},
J. Indian Math. Soc. {\bf 20}(1956), 47-87.

\item{ [56]} R. Sherman Lehman, {\it On the difference $\pi(x) - {\rm
li}\,x$},
Acta Arith, {\bf 11}(1966), 397-410.

\item{ [57]} C.L. Siegel, {\it \"Uber Riemanns {Nachla\ss}  zur analytischen
Zahlentheorie}, Quell. Stud. Gesch. Mat. Astr.
Physik {\bf 2}(1932), 45-80 (also in
{\it Gesammelte Abhandlungen, Band I}, Springer Verlag, Berlin etc.,
1966,  275-310).

\item{ [58]} R. Spira, {\it Some zeros of the Titchmarsh counterexample},
Math. Comp. {\bf 63}(1994), 747-748.

\item{ [59]} K.-M. Tsang, {\it Some $\Omega$-theorems for the Riemann
zeta-function}, Acta Arith. {\bf 46}(1986), 369-395.

\item{ [60]} K.-M. Tsang, {\it The large values of the Riemann zeta-function},
Mathematika {\bf40}(1993), 203-214.

\item{ [61]} E.C. Titchmarsh, {\it The theory of the Riemann zeta-function}
(2nd ed.), Clarendon Press, Oxford, 1986.

\bigskip
\centerline{\bf References added in November 2003:}
\bigskip

\item{ [62]} A. Ivi\'c, {\it On some results concerning the
Riemann Hypothesis}, in
``Analytic Number Theory" (Kyoto,
1996) ed. Y. Motohashi, LMS LNS {\bf247}, Cambridge University Press, 
Cambridge, 1997, 139-167.

\item{ [63]} Y. Motohashi, Spectral  theory of the Riemann zeta-function,
Cambridge University Press, 1997.

\item{ [64]} A. Ivi\'c, {\it On the multiplicity of  zeros of the
zeta-function},
Bulletin CXVIII de l'Acad\'e\-mie Serbe des Sciences et des Arts - 1999, 
Classe des Sciences math\'ematiques et naturelles, 
Sciences math\'ematiques No. {\bf24}, 119-131.

\item {[65]} E. Bombieri,  {\it Problems of the Millenium: the
Riemann Hypothesis},
{\tt http//www.ams.org/claymath/\break prize{\_}problems/riemann.pdf}, 
Amer. Math. Soc., Providence, R.I., 2000, 12pp. 

\item{[66]} J.B. Conrey, {\it $L$-functions and random matrices},
in ``Mathematics Unlimited" (Part I), B. Engquist and W. Schmid eds.,
 Springer, 2001, 331-352.

\item {[67]} J.B. Conrey, D.W. Farmer, J.P. Keating, M.O. Rubinstein
and N.C. Snaith, {\it Integral moments of $L$-functions}, 2003,  58pp,
arXiv:math.NT/0206018, {\tt http://front.math.ucdavis.edu/mat.NT/0206018}.

%\vskip4cm
%\bigskip\bigskip
\vfill

{\cc
\no Aleksandar Ivi\'c

\no Katedra Matematike RGF-a

\no Universitet u Beogradu

\no  \DJ u\v sina 7, 11000 Beograd

\no Serbia (Yugoslavia)}

\no
{\ee aivic@matf.bg.ac.yu,\enskip
aivic@rgf.bg.ac.yu}

\bigskip

\bye